\documentclass[12pt]{amsart}
\font\emailfont=cmtt10

\usepackage{amsmath,amssymb,amsthm,amsfonts,amscd,flafter,epsf,epsfig}

\newcommand\commentable[1]{#1}

\newcommand{\rk}{\mathrm{rk}}

\newtheorem{theorem}{Theorem}[section]
\newtheorem{prop}[theorem]{Proposition}

\newtheorem{lemma}[theorem]{Lemma}

\newtheorem{defn}[theorem]{Definition}

\def\endproof{\relax\ifmmode\expandafter\endproofmath\else
  \unskip\nobreak\hfil\penalty50\hskip.75em\hbox{}\nobreak\hfil\bull
  {\parfillskip=0pt \finalhyphendemerits=0 \bigbreak}\fi}
\def\endproofmath$${\eqno\bull$$\bigbreak}
\def\bull{\vbox{\hrule\hbox{\vrule\kern3pt\vbox{\kern6pt}\kern3pt\vrule}\hrule}}

\newcommand{\Z}{\mathbb{Z}}

\newcommand{\OneHalf}{\frac{1}{2}}

\newcommand{\cm}{\cdot}

\newcommand{\ModSWfour}{\mathcal{M}}
\newcommand{\ModFlow}{\ModSWfour}

\newcommand{\SpinC}{{\mathrm{Spin}}^c}

\newcommand\abuts\Rightarrow
\newcommand\Sym{\mathrm{Sym}}

\newcommand\Cr{\mathrm{Cr}}
\newcommand\HFKm{{\mathrm {HFK}}^-}
\newcommand\HFKa{\widehat{\mathrm{HFK}}}
\newcommand\AlexT{\widetilde\Delta}
\newcommand\ws{\mathbf w}
\newcommand\zs{\mathbf z}

\newcommand\Da{\widetilde{\partial}}
\newcommand\Db{\partial}
\newcommand\HFa{\widetilde{\mathrm{HFS}}}
\newcommand\CFa{\widetilde{\mathrm{CFS}}}
\newcommand\HFb{\mathrm{HFS}}
\newcommand\CFb{\mathrm{CFS}}

\newcommand\x{\mathbf x}
\newcommand\w{\mathbf w}
\newcommand\z{\mathbf z}

\newcommand\y{\mathbf y}

\newcommand\ModSphere{\ModFlow\left({\mathbb S}\longrightarrow 
\Sym^{g-1}(\Sigma_{1})\times \Sym^2(\Sigma_{2})\right)}
\newcommand\ModSpheres\ModSphere

\newcommand\gr{\mathrm{gr}}
\newcommand\Mas{\mu}
\newcommand\UnparModSp{\widehat \ModSp}
\newcommand\UnparModFlow\UnparModSp
\newcommand\Mod\ModSp

\newcommand{\cald}{{\mathcal D}}

\newcommand\ModMaps{\mathcal M}
\newcommand\ModSp\ModMaps

\newcommand\Ta{{\mathbb T}_{\alpha}}
\newcommand\Tb{{\mathbb T}_{\beta}}
\newcommand\Tc{{\mathbb T}_{\gamma}}

\newcommand\alphas{\mbox{\boldmath$\alpha$}}

\newcommand\betas{\mbox{\boldmath$\beta$}}
\newcommand\gammas{\mbox{\boldmath$\gamma$}}

\newcommand\Field{\mathbb F}

\newcommand\Dual{\mathcal D}
\newcommand\Duality\Dual

\title[{Floer homology and singular knots}]
{Floer homology and singular knots}

\author[Peter Ozsv{\'a}th]{Peter Ozsv\'ath}
\address{Department of
Mathematics, Columbia University, New York 10027\newline
\indent{\emailfont{petero@math.columbia.edu}}}
\thanks{PSO was supported by NSF grant number DMS 0234311}

\author[Andr{\'a}s Stipsicz]{Andr{\'a}s Stipsicz}
\address{R\'enyi Institute of Mathematics, Hungarian Academy of Sciences,
H-1053 Budapest, Re\'altanoda utca 13--15, HUNGARY\newline
\indent{\emailfont{stipsicz@math-inst.hu}}}
\thanks{AS was supported by OTKA T49449}

\author[Zolt{\'a}n Szab{\'o}]{Zolt{\'a}n Szab{\'o}} 
\address{Department of  Mathematics, Princeton University, 
Princeton, New Jersey 08544 \newline
\indent{\emailfont{szabo@math.princeton.edu}}} 
\thanks{ZSz was supported by NSF grant number DMS 0107792} 
\thanks{AS and ZSz were partially supported by
  the EU Marie Curie TOK program BudAlgGeo}

\renewcommand{\a}{\mathbf a}
\renewcommand{\b}{\mathbf b}
\newcommand{\D}{\mathbf D}

\begin{document}

\begin{abstract}  
  In this paper we define and investigate variants of the link Floer
  homology introduced by the first and third authors.  More precisely,
  we define Floer homology theories for oriented, singular knots in
  $S^3$ and show that one of these theories can be calculated
  combinatorially for planar singular knots.
\end{abstract} 

\maketitle
\section{Introduction}

Singular knots arise naturally in the skein theory of ordinary knots,
see for example~\cite{moy}. Moreover, they play a crucial role in the
Khovanov-Rozansky categorification of the HOMFLY-PT
polynomial~\cite{KhovanovRozansky}. In this paper we define and
investigate Floer homology theories for oriented, singular knots in
$S^3$, whose Euler characteristic is the Alexander polynomial of the
knot.  The definitions of these theories are generalizations of knot
and link Floer homology,
\cite{Knots}, \cite{RasmussenThesis}, \cite{Links}.  

\begin{figure}
\mbox{\vbox{\epsfbox{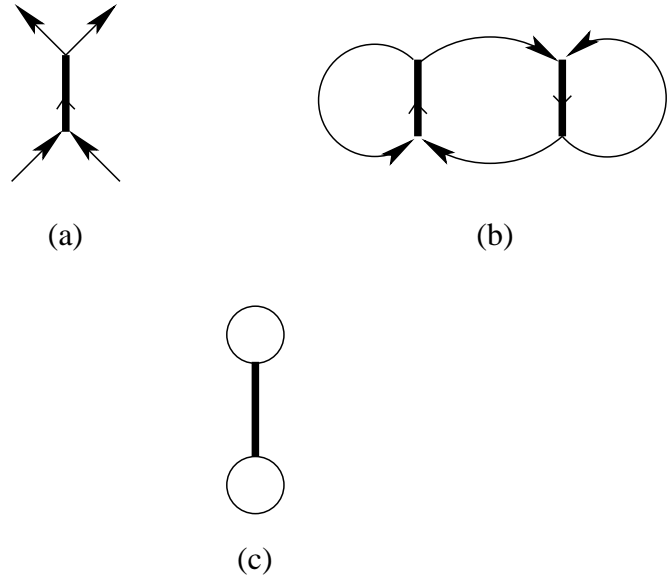}}}
\caption{\label{fig:ThickGraph} 
{\bf{Singular knots.}} (a) shows a thick edge, (b) illustrates the orientation
  convention, while (c) is an example of a singular knot which cannot be
  oriented.}
\end{figure}

Informally, one can think of a singular knot as an ordinary knot with a finite
set of double points. However, for the purpose of skein theory, it is important
to endow these objects with some extra structure.  To this end an {\em
  abstract singular knot} is a connected, trivalent graph with a distinguished
set of edges, called {\em thick edges} (and the remaining edges are called
{\em thin edges}), which satisfies an additional hypothesis: at each vertex,
we require two of the edges to be thin, and the third to be thick.  An {\em
  oriented abstract singular knot} is obtained by orienting all the edges in
such a manner that at each vertex, both thin edges are oriented the same way,
while the thick edge is oriented oppositely, see
Figure~\ref{fig:ThickGraph}(a).  For an example of an oriented singular knot
see Figure~\ref{fig:ThickGraph}(b).  Of course, not all abstract singular
knots can be oriented in this way; an example for a nonorientable abstract
singular knot is shown in Figure~\ref{fig:ThickGraph}(c).  A {\em singular
  knot}, then, is a PL embedding of an oriented abstract singular knot in
$S^3$. By contracting all the thick edges to points, we obtain an oriented
knot with a finite set of double points.

\begin{figure}
\mbox{\vbox{\epsfbox{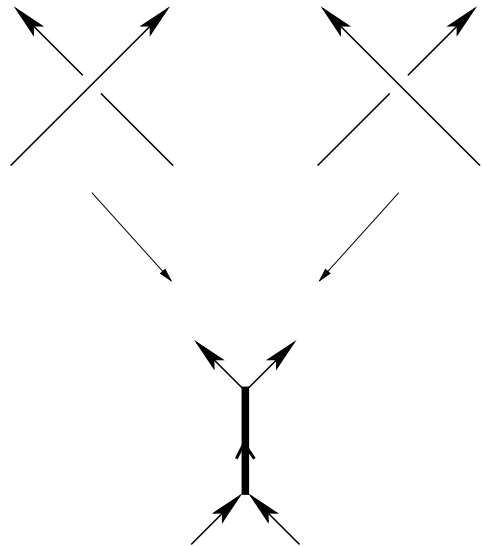}}}
\caption{\label{fig:Singularization} 
{\bf{Singularizing knot projections.}} The diagram shows how to replace a 
  positive/negative crossing by a thick edge in the projection.}
\end{figure}

As it is customary in knot theory, we will mainly study singular knots through
their generic projections. In fact, we consider projections which have
crossings only among thin edges.  Singular knots arise naturally from
projections of oriented knots as follows. Given a projection of an ordinary
knot, we replace some of its crossings by inserting a thick edge, in such a
manner that the orientations of the thin edges are preserved. This is
illustrated in Figure~\ref{fig:Singularization}.

Let $K$ be an oriented, singular knot.  The
Alexander polynomial $\Delta_K(T)$ of $K$ can be defined using
the extension of Alexander polynomials of ordinary links via
skein theory, as it is discussed in  Section~\ref{sec:Alex}.

In this paper we construct a bigraded homology theory associated to
$K$ whose graded Euler characteristic is determined by $\Delta_K(T)$.
In fact, there are two variants of this homology theory, $\HFa(K)$ and
$\HFb(K)$, and each splits as
\begin{eqnarray*}
\HFa(K)=\bigoplus_{s,d} \HFa_d(K,s)
&{\text{and}}&
\HFb(K)=\bigoplus_{s,d} \HFb_d(K,s)
\end{eqnarray*} 
where $d$ (corresponding to the homological, or Maslov grading) is an
integer, and $s$ (the Alexander grading) is either an integer or a
half-integer, depending on the number of thick edges of the singular
knot.  The group $\HFb(K)$ is a module over the polynomial ring
$\Field [U_1,\ldots , U_{\ell }]$, where $\ell +1$ denotes the number
of thick edges of $K$ ($\ell \geq 0$) and $\Field =\Z/2\Z$. 
In fact,
$\HFb(K)$ is gotten as the homology of a chain complex $\CFb(K)$ over
this ring. The group $\HFa(K)$ is obtained by taking the homology of
the complex $\CFa(K)$ gotten by specializing all $U_i=0$ in
$\CFb(K)$. Multiplication by $U_i$ drops Alexander grading by one, and
it drops Maslov grading by two. 

Note that in the definition of $\HFb(K)$, the number of variables is
one less than the number of thick edges; our construction treats one
of these edges specially. Thus, these Floer homology groups depend on
this choice of a distinguished thick edge on the singular knot;
i.e. the resulting group is an invariant of the singular knot equipped
with a thick edge. The groups $\HFa(K)$ use all thick edges more symmetrically,
and as such depends only on the singular link. See Theorem~\ref{thm:WellDefined}
for the invariance statement.

The relationship between these invariants
and the Alexander polynomial is spelled out in the following:

\begin{theorem}
        \label{thm:Euler}
        For the Floer homology theory $\HFa$ we have
        \begin{equation}
          \label{eq:EulerHFa}
          \sum_{s}\chi(\HFa_*(K,s))\cm T^s=
          (1-T)^{\ell}\cm\Delta_K(T),
        \end{equation}
        where
        $\ell +1$ denotes the number of thick edges in $K$. 
        For $\HFb$ we have     
        \begin{equation}
          \label{eq:EulerHFb}
          \sum_{s}\chi(\HFb_*(K,s))\cm T^s=\Delta_K(T).
        \end{equation}
\end{theorem}

A singular knot projection gives rise to a natural Heegaard diagram.
Using such a diagram, and generalizing the notion of Kauffman
states to singular knots (see Section~\ref{sec:States} for details),
we get a natural interpretation of a set of generators for the Floer
homology groups in terms of these extended Kauffman states, similarly
to \cite[Theorem~1.2]{Alternating}. 

\begin{theorem}  \label{thm:States}
  Let $K$ be a singular knot, and fix a decorated projection $P$ of
  $K$.  Consider the vector space $C(K)$ over $\Field$ generated by all generalized
  Kauffman states of $P$.  There is a differential on $C(K)$
  whose homology calculates  $\HFb(K)$.  
\end{theorem}

A more precise statement is given in Section~\ref{sec:States}, where
we also give an explicit formula for the bigrading associated to each
generalized Kauffman state.

We turn our attention to the Floer homologies of singular knots in
the case where the singular knot $K$ is {\em planar}, that is, $K$
admits an injective projection to the plane.

\begin{theorem}
  \label{thm:Planar} Suppose that $K$ is a planar singular knot. Then the
  group $\HFb(K)$ is determined by the Alexander polynomial
  $\Delta_K(T)$; indeed, the homology $\HFb_d(K,s)$ is supported
  on the line $2s=d$.
\end{theorem}

It is more challenging to calculate $\HFa(K)$. In particular, we exhibit an
example which shows that $\HFa(K,s)$ can be non-trivial even when its
Euler characteristic vanishes.

The present paper is organized as follows. In Section~\ref{sec:Defs}
we give the definition of $\HFa$ and $\HFb$. In Section~\ref{sec:Alex}
we briefly discuss Alexander polynomials of singular knots.  In
Section~\ref{sec:States} we give the state sum model for the Floer
homology theory $\HFb$ of singular knots, leading us to the proofs of
the theorems announced above. In Section~\ref{sec:Planar} we show an
alternative Heegaard diagram for planar singular knots and do some
calculations in some special cases.  

Most of the material described in this paper was discovered in 2004,
while the second author was visiting the Institute for Advanced
Study. The invariants for singular knots described here (and some
suitable modification of them) form the basis for a {\em cube of
resolutions} description of knot Floer homology, described
in~\cite{Resolutions}, which gives a construction of knot Floer
homology quite similar in character to the Khovanov-Rozansky
categorification of the HOMFLY-PT polynomial~\cite{KhovanovRozansky}.
A different construction of knot Floer homology for singular links has
appeared very recently in~\cite{Audoux}, cf. also \cite{Nadya}.

\subsection{Acknowledgements} The authors wish to thank the referee for
a careful reading and extensive suggestions.

\section{Definition of the Floer homology groups}
\label{sec:Defs}

\subsection{Heegaard diagrams}
\label{sec:HeegaardDiag}

Consider an oriented surface $\Sigma$ of genus $g$, endowed with two
$(g+\ell)$-tuples of mutually pairwise disjoint, embedded circles
$\alphas=\alpha_1,\ldots ,\alpha_{g+\ell}$ and
$\betas=\beta_1,\ldots  ,\beta_{g+\ell}$ (where $\ell$ is a non-negative
integer), 
chosen so that $\alpha_1, \ldots  ,\alpha_{g+\ell}$
and $\beta_1, \ldots ,\beta_{g+\ell}$ determine handlebodies $U_\alpha$ and
$U_\beta$ with $\partial U_{\alpha}= -\partial U_{\beta}=\Sigma$. If
$U_\alpha\cup_{\Sigma} U_\beta\cong S^3$, we call
$(\Sigma,\alphas,\betas)$ a {\em balanced Heegaard diagram for $S^3$}.
A self-indexing Morse function on $S^3$ with $\ell+1$ maxima and
minima determines a balanced Heegaard diagram for $S^3$.

A point $p\in \Sigma-\alpha_1- \ldots -\alpha_{g+\ell}-\beta_1- \ldots
-\beta_{g+\ell}$ determines a gradient flow which connects some index zero and
index three critical point.  Indeed, given a collection of points
$\ws=w_1,\ldots ,w_a$ and $\zs=z_1,\ldots ,z_b$, each in $\Sigma-\alpha_1-
\ldots -\alpha_{g+\ell}-\beta_1- \ldots -\beta_{g+\ell}$, we obtain an
oriented graph $\Gamma_{\ws,\zs}$, which is oriented compatibly with the flow
along the edges corresponding to $w_i$, and oppositely for the edges
corresponding to $z_j$.  We suppose moreover that $\Gamma_{\ws,\zs}$ is a
connected graph with vertices of degree three. This condition implies that
each component $\{B_i\}_{i=1}^{\ell+1}$ of $\Sigma-\beta_1- \ldots
-\beta_{g+\ell}$ has two points of type $\zs$ (denoted by $z^1_i$ and
$z^2_i$), and one point of type $\ws$ (which will be denoted by $w_i$).
A similar statement holds for the components $\{ A_j \} _{j=1} ^{\ell +1}$ of
$\Sigma-\alpha_1- \ldots -\alpha_{g+\ell}$.  In this way we obtain an oriented
singular knot in $S^3$ by declaring the edges passing through $w_{1},\ldots
,w_{\ell+1}$ to be thick edges.  Conversely, it is not hard to see that any
singular knot $K\subset S^3$ arises in this way, that is, for any $K$ there is
a balanced Heegaard diagram compatible with it. An example of such a Heegaard
diagram will be given in the next subsection.

\subsection{Decorated projections and Heegaard diagrams}
\label{sec:HeegaardProjection}

\begin{figure}
\mbox{\vbox{\epsfbox{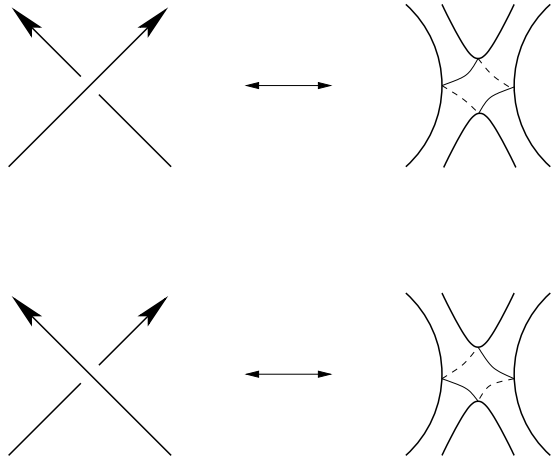}}}
\caption{\label{fig:Crossing}
{\bf{Heegaard diagram at a crossing.}} The position of the $\alpha$ circle
(depending on the sign of the crossing) is shown.}
\end{figure}

\begin{figure}
\mbox{\vbox{\epsfbox{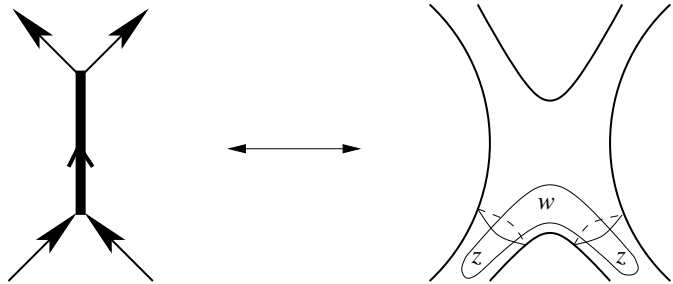}}}
\caption{\label{fig:HeegSingPoints} {\bf{Heegaard diagram at a thick edge.}}
  The diagram describes the $\alpha$-cirles and the internal $\beta$-circle
  corresponding to a thick edge, together with the base points $w$ and two
  $z$'s.}
\end{figure}

Fix a generic projection $P$ of an oriented singular knot $K$.  We
also fix a generic initial point $Q$ on the projection. We call this
data a {\em decorated knot projection} for $K$.  Given this data, we
describe an associated Heegaard diagram, generalizing the one given in
\cite{Alternating}.  (Note that we have switched here the roles of the
$\alpha$ and $\beta$ circles from the conventions
from~\cite{Alternating}.)  The handlebody $U_\alpha$ is a regular
neighborhood of the projection of $K$, and $U_\beta$ is its complement
in $S^3$.  Let $X$ and $Y$ be the regions in the projection which
contain the distinguished point $Q$.  For each region in the
complement of the projection other than $X$, we choose a corresponding
$\beta$-circle given by the contour of the region.  For each crossing
in the projection we choose a circle $\alpha_j$ supported in a
neighborhood of the crossing as it is shown in
Figure~\ref{fig:Crossing}.
For each thick edge, we choose a pair of $\alpha$-circles $\alpha_{j}$
and $\alpha_{j+1}$ which are meridians for the two incoming arcs, and
also an additional circle $\beta_i$, which meets only $\alpha_j$ and
$\alpha_{j+1}$ in two points apiece. We call $\beta_i$ an {\em
  internal $\beta$-circle} for the thick edge.  This is illustrated in
Figure~\ref{fig:HeegSingPoints}.  The diagram also shows how to place
the base points $w_i$ and $z_i^1, z_i^2$ near the thick edge.  To
complete the construction, finally we omit one of the internal
$\beta$-circles. We will take the convention that $w_{\ell+1}$ is the
distinguished thick edge. Moreover, in Section~\ref{sec:States}, we
will find it useful to choose a diagram where the initial point $Q$
lies on a thin edge which points into the distinguished thick edge,
and that the internal $\beta$-circle for this thick edge is omitted.
It is easy to see that the resulting balanced Heegaard diagram is
compatible with the given singular knot.

\subsection{Definition of  $\HFa(K)$ and $\HFb(K)$}
Suppose now that $K\subset S^3$ is a singular knot and
$(\Sigma,\alphas,\betas,\w, \z)$ is a balanced Heegaard diagram
compatible with $K$.  As usual (compare~\cite{HolDisk}, \cite{Links}),
we consider the $(g+\ell)$-fold symmetric power $Sym ^{g+\ell }(\Sigma )$
of the genus-$g$ surface $\Sigma$, equipped with the tori
\begin{eqnarray*}
\Ta=\alpha_1\times\ldots \times\alpha_{g+\ell}
&{\text{and}}&
\Tb=\beta_1\times \ldots \times\beta_{g+\ell}
\end{eqnarray*}
in $Sym ^{g+\ell }(\Sigma )$.  Let $\pi _2 (\x , \y )$ denote the
set of homology classes of topological disks connecting the intersection
points $\x, \y \in \Ta \cap \Tb$. For $p\in \Sigma-\alpha_1- \ldots
-\alpha_{g+\ell}-\beta_1- \ldots -\beta_{g+\ell}$ and $\phi \in \pi _2
(\x , \y )$ we define $n_p(\phi )$ as the intersection of $\phi$ with
$V_p=\{ p\}\times Sym ^{g+\ell-1}(\Sigma)$. We decompose the
complement of the $\alpha_i$ and $\beta_j$ into its components
$$\Sigma-\alpha_1- \ldots -\alpha_{g+\ell}-\beta_1- \ldots
-\beta_{g+\ell}=\coprod_k C_k.$$ Given $\phi\in\pi_2(\x,\y)$, the {\em domain}
$\cald(\phi)$ associated to $\phi$ 
is the formal linear combination
$\cald(\phi) = \sum _j n_j C_j$, where here $n_j=n_{c_j}(\phi)$ for any choice
of $c_j\in C_j$. The homology class of $\phi$ is uniquely determined by its
associated domain. We say that the domain $\cald$ is \emph{nonnegative} (writing
$\cald \geq 0$) if $\cald = \sum _j n_j C_j$ with $n_j \geq 0$ for all
$j$.
  
Two functions
$$A'\colon (\Ta\cap\Tb)\times (\Ta\cap\Tb) \longrightarrow \Z$$
and 
$$N'\colon (\Ta\cap\Tb)\times (\Ta\cap\Tb) \longrightarrow \Z$$
can be defined by the formulas
\begin{equation}
  \label{eq:AlexanderDifference}
A'(\x,\y)=
\sum_{i=1}^{\ell+1} \left(2n_{w_i}(\phi)-n_{z^1_i}(\phi)-n_{z^2_i}(\phi)
\right)
\end{equation}
and 
$$N'(\x,\y)=\Mas(\phi)-2\sum_{i=1}^{\ell+1} n_{w_i}(\phi),$$ where
$\phi \in \pi _2 (\x, \y ) $ is arbitrary, and $\mu (\phi )$ is the
Maslov index of $\phi \in \pi _2 (\x, \y )$.

\begin{lemma}
\label{lemma:MaslovSpinC}
The functions $A'$ and $N'$ are well-defined, i.e. they are independent
of the choice of $\phi\in\pi_2(\x,\y)$.
Moreover, given $\x,\y,\w\in\Ta\cap\Tb$, we have that
$A'(\x,\y)+A'(\y,\w)=A'(\x,\w)$ and
$N'(\x,\y)+N'(\y,\w)=N'(\x,\w)$.
\end{lemma}

\begin{proof}
  Consider $\phi$ and $\phi'\in \pi_2(\x,\y)$.  The difference
  $D=\cald(\phi)-\cald(\phi')$ of their domains specifies a two-chain whose
  boundary is a sum of $\alpha$- and $\beta$-circles.  Since $H_1(\Sigma)$ is
  spanned by the images of the two $g$-dimensional subspaces spanned by the
  $\{[\alpha_i]\}_{i=1}^{g+\ell}$ and $\{[\beta_i]\}_{i=1}^{g+\ell}$, it
  follows that any such domain can be written as a sum of some $D_1$ and
  $D_2$, where
\begin{eqnarray*}
 D_1 = \sum_{i=1}^{\ell+1} a_i\cm A_i
&{\text{and}}&
 D_2 = \sum_{i=1}^{\ell+1} b_i\cm  B_i,
\end{eqnarray*}
where $\{A_i\}_{i=1}^{\ell+1}$ resp.  $\{B_i\}_{i=1}^{\ell+1}$ are the
components of $\Sigma-\alpha_1-\ldots -\alpha_{g+\ell}$ resp.
$\Sigma-\beta_1-\ldots -\beta_{g+\ell}$, and $a_i,b_i\in\Z$.  Now, our
combinatorial condition on the Heegaard diagram ensures that each component
$A_i$ or $B_j$ contains two points of $\z$ and one point of $\w$
(corresponding to a thick edge). Since these latter points contribute to the
sum in $A'$ with multiplicity two, it follows readily that
\begin{eqnarray*}
\lefteqn{
\sum_{i=1}^{\ell+1} \left(2n_{w_i}(\phi)-n_{z^1_i}(\phi)-n_{z^2_i}(\phi)\right) }\\ 
&=&
\sum_{i=1}^{\ell+1} \left(2n_{w_i}(\phi')-n_{z^1_i}(\phi')-n_{z^2_i}(\phi')\right);
\end{eqnarray*}
i.e. $A'$ is well-defined.
Well-definedness of $N'$ follows from standard properties of the
Maslov index, see~\cite[Proposition~4.1]{Links}. 
Additivity of these quantities is straightforward.
\end{proof}

It follows from Lemma~\ref{lemma:MaslovSpinC} that there are functions
$A\colon \Ta\cap\Tb \longrightarrow \Z$ and $N\colon \Ta\cap\Tb
\longrightarrow \Z$, uniquely characterized up to an overall translation by
the property that $A(\x)-A(\y)=A'(\x,\y)$ and $N(\x)-N(\y)=N'(\x,\y)$.  It
will be useful to have the following

\begin{lemma}
  \label{lemma:UniqueHoClass}
  Given $\x\in\Ta\cap\Tb$, there is a unique $\phi\in\pi_2(\x,\x)$
  with $n_{\ws}(\phi)=n_{\zs}(\phi)=0$.
\end{lemma}
\begin{proof}
  Notice first that the trivial class $\phi _0\in \pi _2 (\x, \x)$ has the
  required property. Any other class in $\pi _2 (\x, \x)$ differs from $\phi
  _0$ by adding $D_1$ and $D_2$ as above. Now the condition on $n_{\ws}$ and
  $n _{\zs}$ implies that for any $a_i$ in $D_1$ and $b_j$ in $D_2$ with the
  property that $A_i \cap B_j$ contains some $z_i$ or $w_i$ we have that
  $b_j=-a_i$.  Using this repeatedly, together with the fact that
  $\Gamma_{\ws,\zs}$ is connected, we see that $D_1+D_2$ is zero,
  concluding the proof.
\end{proof}

We define the complex $\CFb(K)$ to be the chain complex freely
generated over $\Field [U_1,\ldots ,U_{\ell}]$ (with $\Field = \Z
/2\Z$) by the intersection points $\Ta\cap\Tb$, endowed with the
differential
\begin{equation}
  \label{eq:DefDb}
\Db (\x)
=\sum_{\y\in\Ta\cap\Tb}\sum_{\{\phi\in\pi_2(\x,\y)\big|\Mas(\phi)=1,
  n_{\w}(\phi)=0\}} \#\UnparModFlow(\phi) \cdot
\left(\prod_{i=1}^{\ell} U_i^{(n_{z^1_{i}}(\phi)
    +n_{z^2_{i}}(\phi))} \right) \y,
\end{equation}
where $\#\UnparModFlow(\phi)$ is the mod 2 count of holomorphic
representatives of $\phi \in \pi _2 (\x, \y )$ up to
reparametrization.  

The functions $A$ and $N$ defined above induce gradings on $\CFb(K)$,
the {\em Alexander} and {\em algebraic} gradings, with the convention
that multiplication by $U_i$ drops $A$-grading by one, and preserves
$N$.

For the case of non-singular knots, knot Floer homology is bigraded,
with Alexander and Maslov degrees. The Maslov degree of~\cite{Knots}
behaves differently from the grading $N$ introduced here:
multiplication by $U_i$ drops Maslov grading, while it preserves $N$.
In fact, we could define the following {\em Maslov grading} for the
case of singular links by taking
$$
M'(\x,\y)=\Mas(\phi)+2\cm\sum_{i=1}^{\ell+1}\left(n_{w_i}(\phi)-n_{z^1_i}(\phi)-n_{z^2_i}(\phi)\right).
$$
This is analogous to the Maslov grading for knots, and in
particular, with respect to the induced grading on $\CFb(K)$,
multiplication by any $U_i$ drops this grading by two.  Note that we
have the following easily checked relationship between the three gradings:
$$N'=M'-2A'.$$

The complex $\CFa(K)$ is the chain complex we get by setting all
$U_i=0$ ($i=1, \ldots , \ell$) in $\CFb(K)$.  Equivalently, it is the
free Abelian group generated over $\Field$ by $\Ta\cap\Tb$, endowed with
the differential
$$\Da (\x)
=\sum_{\y\in\Ta\cap\Tb}\sum_{\{\phi\in\pi_2(\x,\y)\big|\Mas(\phi)=1,
n_{\ws}(\phi)=0,
n_{\zs}(\phi)=0\}} \#\UnparModFlow(\phi) \cm \y.$$

\begin{lemma}
  The maps $\Da$ and $\Db$ induce differentials on
  $\CFa$ and $\CFb$ which drop the $N$-grading by one, and preserve
  the $A$-grading.
\end{lemma}

\begin{proof} 
  The sum defining $\Da$ is finite, according to
  Lemma~\ref{lemma:UniqueHoClass}. For $\Db$, we argue that for fixed
  $\x$ and $\y$, there are only finitely many $\phi\in\pi_2(\x,\y)$
  with $n_{\ws}(\phi)=0$ and $\cald(\phi) \geq 0$. To see this,
  observe that the homology class of any $\phi$ with $n_{\ws}(\phi)=0$
  is uniquely determined by $n_{\zs}(\phi)$ (by
  Lemma~\ref{lemma:UniqueHoClass}), and moreover according to
  Lemma~\ref{lemma:MaslovSpinC}, the sum
  $$\sum_{i=1}^{\ell+1} n_{z^1_i}(\phi)+n_{z^2_i}(\phi)$$
  is independent of the choice of $\phi$ (since it is $A(\y)-A(\x)$).

Since the coefficients in our theories are chosen in $\Field = \Z
/2\Z$, the arguments of \cite{Links} apply. Now the proof that
$\Da^2=0=\Db^2$ is an adaptation of \cite[Lemma~4.3]{Links}.
\end{proof}

Note that, we have so far defined $N$ and $A$ only up to an overall
additive constant.  The indeterminacy in $N$ can be easily removed
with the following observation.  Suppose that we set all the $U_i=1$
in the definition of $\CFb$ (equivalently, consider the differential
which counts holomorphic disks with $n_{\ws}(\phi)=0$, and ignore the
reference points $\zs$). In this way we obtain a chain complex whose
homology is (up to an overall degree shift) isomorphic to
$H_*(T^\ell)$,
see~\cite[Theorem~\ref{Links:thm:InvarianceHFaa}]{Links}.  Note that
by Lemma~\ref{lemma:UniqueHoClass}, it follows that this diagram is
admissible, so that the theorem applies.
We choose $N$ so that the homology group is isomorphic to
$H_*(T^\ell)$, with a shift in the grading making the top-most
homology supported in degree zero.

The indeterminacy of $A$ can be removed in an invariant manner using
relative $\SpinC$ structures (compare~\cite{TuraevTorsion}), though
this is a slightly awkward approach for the present applications.
Rather, we will leave the $A$-grading well-defined only up to an
overall shift, which could be removed with the help of the state sum
formulae of Section~\ref{sec:States}. We now write
\begin{eqnarray*}
\CFa(K)=\bigoplus_{d,s}\CFa_d(K,s) &{\text{and}}&
\CFb(K)=\bigoplus_{d,s}\CFb_d(K,s)
\end{eqnarray*}
where here $\CFa_d(K,s)$ and $\CFb_d(K,s)$ are generated by elements with
$M=d\in \Z$ and $A=s$. (Presently, the grading by $s$ is well-defined only up
to an additive constant.)

\begin{theorem}
  \label{thm:WellDefined}
  Consider the homology groups
  \begin{eqnarray*}
    \HFa(K)=\bigoplus_{d,s} \HFa_d(K,s)
    &{\text{and}}&
    \HFb(K)=\bigoplus_{d,s} \HFb_d(K,s)
  \end{eqnarray*}
  of the chain complexes $(\CFa (K), \Da)$ and $(\CFb (K), \Db)$, where
  here $d\in\Z$ is an absolute integral grading and $s$ (which is in
  $\Z$ or $\Z +\frac{1}{2}$, depending on the parity of the number of
  thick edges) is a relative integral grading.
  Then, $\HFb(K)$ is an invariant of the singular knot $K$, together
  with the distinguished thick edge; while
  $\HFa(K)$ is an invariant of the singular knot.
\end{theorem}

\begin{proof} 
  By Morse theory it is known that any two pointed Heegaard diagrams
  for a fixed singular knot can be connected by a sequence of
  (pointed) handleslides, isotopies, stabilizations and birth/death of
  index 0/1 (and 2/3) cancelling handles. Using the proof of
  \cite[Lemma~1.1]{SharleThomps} we can avoid the appearance of
  cancelling 0/1 and 2/3 handles.  For the remaining moves the proof
  of the invariance of the knot Floer homology groups follows the same
  lines as the proof of \cite[Theorem~4.7]{Links}.  
\end{proof}

\subsection{Other constructions}

There are other variants of Floer homology for singular knots, which
generalize the constructions we describe here. For example, we could
introduce a variable $U_{\ell+1}$ for the final thick edge. The
corresponding theory is the analogue for singular knots of $\HFKm$ for
knots (while the construction $\HFb$ here plays the role of $\HFKa$).

Indeed, this can be further generalized as follows: consider the chain
complex $C'$ which is freely generated by $\Ta\cap\Tb$ over the ring
$\Field [U_1,...,U_m]$, where $m$ denotes the number of {\em thin edges} in our
diagram, and endowed with the differential:
$$\partial' (\x)
=\sum_{\y\in\Ta\cap\Tb}\sum_{\{\phi\in\pi_2(\x,\y)\big|\Mas(\phi)=1,
  n_{\w}(\phi)=0\}} \#\UnparModFlow(\phi) \cdot
\left(\prod_{i=1}^{m} U_i^{n_{z_i}}(\phi)\right)\y.
$$
  
The chain complex $(\CFb,\Db)$ described earlier can be thought of as
the chain complex over the quotient ring where we set $U_{e}=U_{e'}$,
if $e$ and $e'$ are two edges which point into the same vertex; and
also for some edge $e_0$, we set $U_{e_0}=0$. See~\cite{Resolutions} for
related constructions.

\section{Alexander polynomial of singular links}
\label{sec:Alex}

State sum formulas for the Alexander polynomial of
a smooth knot were introduced in~\cite{Kauffman}. 
Our aim in this section is to recall the Alexander
polynomial for singular links, and describe a state sum
description for it, cf. also~\cite{moy}.

\begin{figure}
\mbox{\vbox{\epsfbox{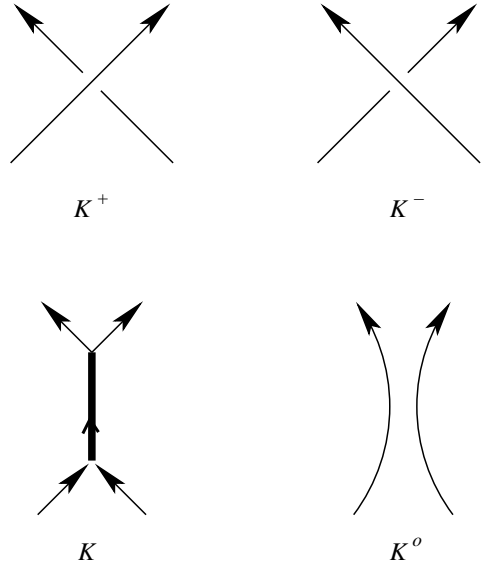}}}
\caption{\label{fig:Resolv} 
{\bf{Resolutions at a thick edge.}}}
\end{figure}

\begin{prop}
        \label{prop:SkeinRelation} There is a unique extension 
        ${\widetilde \Delta}_K(T)$ of the one-variable symmetrized
        Alexander polynomial $\Delta_K(T)$ for non-singular, oriented
        links to singular, oriented links, which is characterized by
        the skein relations
\begin{align*}
        {\AlexT}_{K^+}(T)&=\AlexT_{K} (T)+T^{\OneHalf}\cm {\AlexT}_{K^o}(T)
         \\ 
         \AlexT_{K^-}(T)&=\AlexT_{K}(T)+T^{-\OneHalf}\cm
        \AlexT_{K^o} (T), \\
\end{align*}
        where here $K^o$ denotes the resolution of $K$ at a singular
        point $v$, and $K^+$ and $K^-$ denote the positive and
        negative resolutions at $v$ (cf. Figure~\ref{fig:Resolv}).
\end{prop}

\begin{proof}
  Uniqueness of such an extension follows from simple induction on the
  number of singular points, since $K^{\pm }$ and $K^o$ have fewer
  singular points than $K$. The fact that $\AlexT _K (T)$ is
  well-defined follows from a similar inductive argument together with
  the fact that, because of the ordinary skein relation for the
  Alexander polynomial, the two equations above provide coherent
  results.
\end{proof}

\begin{figure}
\mbox{\vbox{\epsfbox{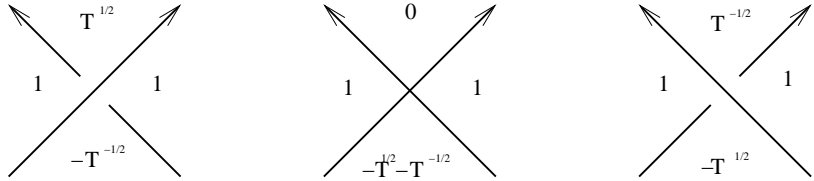}}}
\caption{\label{fig:LocalContrib} 
{\bf{Local contributions.}}  We
illustrate here the local contribution function $Z_v(x)$ of a Kauffman
state at a given crossing $v$. The left and right diagrams illustrate the
case where the crossing is a non-singular point of $K$, while the
middle one illustrates the case where $v$ is a singular point.}
\end{figure}

From now on, the Alexander polynomial ${\widetilde \Delta}_K(T)$ for a
singular knot $K$ will be denoted simply by ${\Delta}_K(T)$.  Now we
turn to the state sum description of $\Delta_K (T)$ for singular
knots.  To this end, first we recall the definition of Kauffman states
for a singular link projection.  Fix a decorated projection $P$ for a
singular link $K$ with initial point $Q$ and contract the thick edges
to points.  Let $\Cr(P)$ and $R(P)$ denote the set of crossings and
the set of regions in the complement of the projection, resp.  Let $X$
and $Y$ be the two regions which contain the edge containing $Q$.  A
{\em Kauffman state} is a bijection $x\colon \Cr(P)\longrightarrow
R(P)-X-Y$ which assigns to each crossing $v\in \Cr(P)$ one of the (up
to) four regions which meet at $v$. Let $x$ be a Kauffman state and
$v$ a crossing in the projection. There is a {\em local contribution}
$Z_v(x)$ defined according to which quadrant $x$ assigns to $v$, and
the type of $v$ (i.e., whether it is singular, positive, or
negative). When $v$ is non-singular, this local contribution takes
values $Z_v(x)\in \{1, \pm T^{\pm \OneHalf}\}$, and when $v$ is
singular, $Z_v(x)\in \{0,1,-T^{\OneHalf}-T^{-\OneHalf}\}$.  The
precise rules are illustrated in Figure~\ref{fig:LocalContrib}.  This
definition can be regarded as an appropriate extension of the local
contributions for Kauffman states of ordinary links.  In the case of
ordinary links, the appropriate sum of the local contributions of
Kauffman states provides the Alexander polynomial;
Proposition~\ref{prop:StateSumFormula} generalizes this fact to
singular links.  We note here that for our later purposes a refinement
of Kauffman states (which we will call \emph{generalized Kauffman
  states}) will be introduced in Section~\ref{sec:States}.

\begin{prop}
        \label{prop:StateSumFormula} 
        Fix a decorated projection $P$ of the oriented link $K$, let $\Cr(P)$
        denote the set of crossings in the projection while $X(P)$ denotes its
        set of Kauffman states. Then the Alexander polynomial $\Delta_K(T)$ is
        calculated using the state sum formula
        \begin{equation}\label{e:sum} \Delta_K(T)= \sum_{x\in X(P)}
        \left(\prod_{v\in \Cr(P)} Z_v(x)\right).  \end{equation}
\end{prop}

\begin{figure}
\mbox{\vbox{\epsfbox{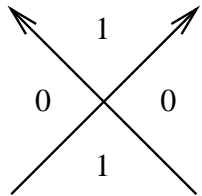}}}
\caption{\label{fig:ResolvedContrib}
{\bf{Local contributions at a resolved point.}}}
\end{figure}

\begin{proof}
        The state sum formula for $K^o$ can be thought of as given by
        a state sum formula for states of $K$, where the local
        contribution at the resolved point contributes as in
        Figure~\ref{fig:ResolvedContrib}.  Since the state sum formula
        satisfies the skein relation of
        Proposition~\ref{prop:SkeinRelation}, and the identity of
        the proposition is known to hold for a non-singular link, the
        formula of \eqref{e:sum} for $\Delta_K(T)$ follows at once.
\end{proof}

 \section{Kauffman states and the chain complex}
\label{sec:States}

Our aim is to give a simple set of generators for a chain complex
whose homology is isomorphic to $\HFb(K)$. These generators can be
described concretely in terms of a suitable generalization of Kauffman
states, which we describe presently.

Fix a decorated projection $P$ for a singular knot $K$ with initial
point $Q$.  We contract all the thick edges, so that they become
crossings of the projection. However, the four quadrants around a
singular point $v$ do not play equal roles. There are two quadrants,
the {\em side quadrants} which correspond to regions which used to
contain the thick edge in their boundary, which we denote $A_v$ and
$C_v$.  There are two remaing quadrants: the {\em top quadrant} $B_v$,
which is pointed towards by the thick edge, and the {\em bottom
  quadrant} $D_v$. We define four quadrants $A_v$, $B_v$, $C_v$, and
$D_v$ at the ordinary crossings analogously.

We define a set, the set of {\em Kauffman corners at $v$}, for each
crossing $v$ of $K$, which depends on whether $v$ is an ordinary
crossing or a contracted thick edge.  If $v$ is an ordinary crossing,
the Kauffman corners are the four corners $A$, $B$, $C$,
or $D$ of the crossing $v$. If $v$ is a singular crossing, the
Kauffman corners correspond to $A$, $C$, $D^+$, and
$D^-$, where both $D^+$ and $D^-$ belong to the bottom corner $D$.

Let $R(P)$ denote the set of regions in the complement of the
projection. Let $X$ and $Y$ be the two regions which contain the edge
containing $Q$. Let $\Cr(P)$ denote the set of crossings of $P$ (i.e.
$\Cr(P)$ consists of the ordinary crossings and also the contracted
thick edges).  A {\em generalized Kauffman state} for a singular knot
is a map $x$ which associates to each crossing $v\in \Cr(P)$ one of
its four allowed Kauffman corners, with the constraint that in each
allowed region in $R(P)$, there is a unique Kauffman corner in the
image of $x$, compare~\cite{Kauffman}.  (Note that this is a very mild
generalization of the notion from Section~\ref{sec:Alex}: there, we
had four states corresponding to four quadrants, except one corner in
the singular case contributed $-T^{\OneHalf}-T^{-\OneHalf}$; for our
present purposes, it is convenient to think of both terms as
corresponding to two different states $D^+$ and $D^-$.)

\begin{figure}
\mbox{\vbox{\epsfbox{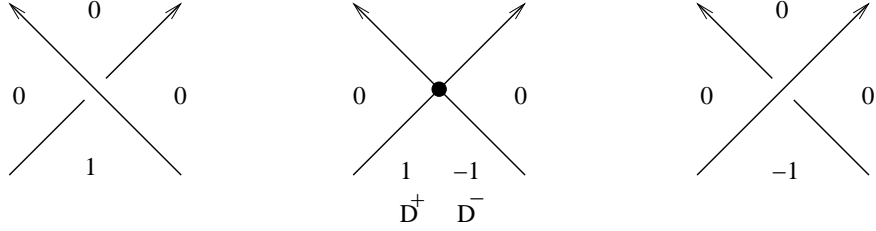}}}
\caption{\label{fig:MasGradings}
{\bf{The local Maslov grading $M_v$ is illustrated at a crossing $v$.}}}
\end{figure}

Each Kauffman corner has a local Maslov grading $M_v$, which is zero on all
Kauffman corners except for $D$ (and $D^+, D^-$ at a singular point), where it
is $\mp 1$ according to the sign of the crossing.  This choice of signs is
illustrated in Figure~\ref{fig:MasGradings}.  The Maslov grading of a
generalized Kauffman state is the sum of local Maslov gradings over each
crossing of $K$. (Notice that at the singular point no value at $B$ is specified, since when
$v$ is a singular point, $B$ is not a Kauffman corner.)

\begin{figure}
\mbox{\vbox{\epsfbox{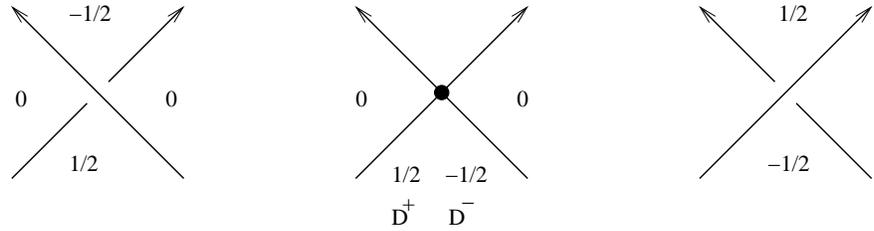}}}
\caption{\label{fig:AlexGradings}
{\bf{The local Alexander grading $S_v$ is illustrated at a crossing $v$.}}}
\end{figure}

Similarly, each Kauffman
corner has a local Alexander grading $S_v$, defined as follows. If $v$ is an
ordinary positive crossing then $S_v(A)=S_v(C)=0$, $S_v(B)=1/2$ and
$S_v(D)=-1/2$; while if $v$ is negative, then $S_v(A)=S_v(C)=0$, $S_v(B)=-1/2$
and $S_v(D)=1/2$; finally, if $v$ is a singular crossing, then
$S_v(A)=S_v(C)=0$ and $S_v(D^-)=-1/2$, $S_v(D^+)=1/2$, see
Figure~\ref{fig:AlexGradings}.  The Alexander grading $S(x)$
of a Kauffman state $x$ is
given as a sum of local Alexander gradings.

\begin{figure}
\mbox{\vbox{\epsfbox{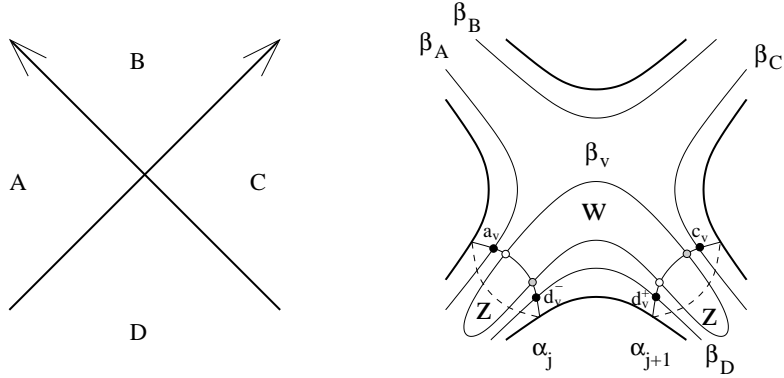}}}
\caption{\label{fig:LocalPicture}
{\bf{The local picture of the Heegaard diagram near a singular point $v$.}}
There are four unlabelled intersection points (the intersection points
of $\beta_v$ with $\alpha_j$ and $\alpha_{j+1}$): two are white and two 
are gray. The white ones correspond to generators with $\epsilon(v)=+1$,
the gray ones to those with $\epsilon(v)=-1$.}
\end{figure}

Fix a decorated projection for a singular knot $K$, with the property
that the distinguished edge $Q$ is just underneath a singular crossing
$p$.  Consider next the associated Heegaard diagram for $K$ described
in Section~\ref{sec:HeegaardProjection}. Leave out the internal
$\beta$-circle $\beta_p$ associated to the singular crossing $p$ which
lies above $Q$. Moreover, $p$ corresponds to the special thick
edge. It is easy to see that there is a Kauffman state associated to
each $\x\in\Ta\cap\Tb$, determined as follows. Suppose that $v$ is a
singular crossing other than $p$, and let $\alpha_j$ and
$\alpha_{j+1}$ be the two circles supported in a neighborhood of the
singular point $v$. Let $\beta_v$ be the corresponding internal
$\beta$ circle, and write $\beta_A$, $\beta_B$, $\beta_C$, and
$\beta_D$ be the four $\beta$-circles corresponding to the four
quadrants meeting at $v$.  Our tuple $\x$ contains exactly one of the
four points
$$\{a_v,c_v,d^+_v,d^-_v\}=(\alpha_j\cup\alpha_{j+1})\cap
(\beta_A\cup\beta_C\cup\beta_D),$$ where $a_v\in \beta_A$, $d^\pm_v\in
\beta_D$, and $c_v\in\beta_C$ (note that the intersection cannot
contain more points since $\beta_j=\beta _v$ meets only $\alpha_j$ and
$\alpha_{j+1}$). We distinguish $d^+_v$ and $d^-_v$ by the following
convention: $d^+_v$ lies on the same $\alpha$-circle as $c_v$.  Now,
our Kauffman state is determined by $x(v)=A_v$, $D^+_v$, $D^-_v$, and
$C_v$, if $\x$ contains $a_v$, $d^+_v$, $d^-_v$, and $c_v$
respectively.  The local picture at such a generic singular crossing
is illustrated in Figure~\ref{fig:LocalPicture}.  At the singular
point $p$, a similar but somewhat simpler picture applies (the
internal $\beta$-circle is missing, as is one of the $\beta$-arcs
corresponding to a distinguished edge). The Kauffman state at an
ordinary crossing is determined similarly, cf.~\cite{Alternating}.

\begin{figure}
\mbox{\vbox{\epsfbox{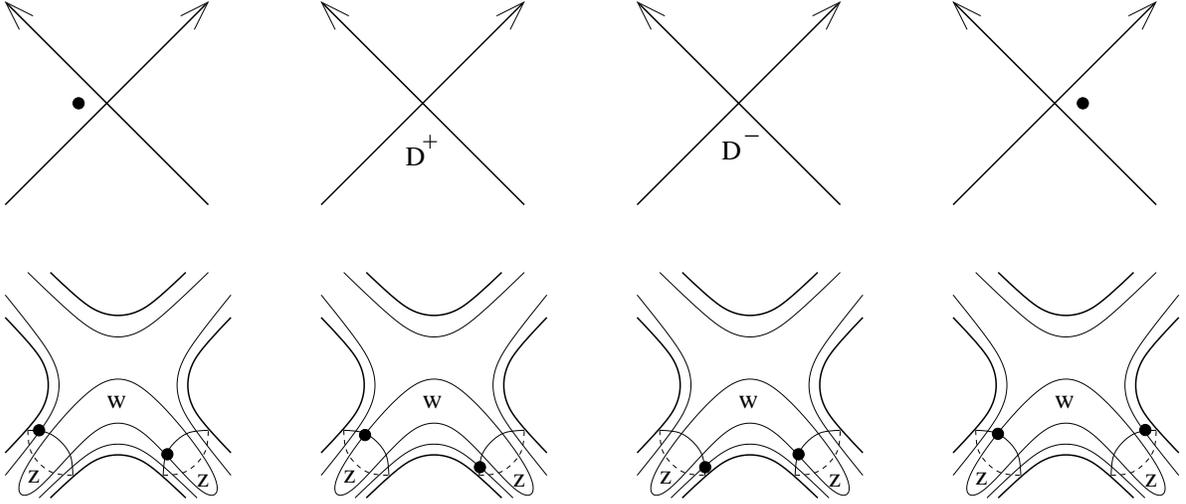}}}
\caption{\label{fig:AlexanderMinimizers} {\bf{Alexander grading
      maximizers.}} For each Kauffman corner illustrated in the top
  row, we have drawn its corresponding Alexander grading maximizing
  intersection point below it.}
\end{figure}

Of course, each Kauffman state is associated to $2^{\ell}$ different
intersection points of $\Ta\cap\Tb$ (where $\ell+1$ denotes the number
of singular points), distinguished by their coordinates on the
internal $\beta$-circles. More specifically, the intersection points
in $\Ta\cap\Tb$ are in one-to-one correspondence with pairs $x$ a
generalized Kauffman state and $\epsilon\colon s(P)\longrightarrow
\{\pm 1\}$, where here $s(P)$ denotes the set of singular points in
the projection near which we have an internal $\beta$-circle (i.e. all
but the distinguished one, $p$). Suppose that $(x,\epsilon)$ and
$(x,\epsilon')$ represent two different intersection points
$\x,\x'\in\Ta\cap\Tb$ with the same underlying Kauffman state, then we
can find a homology class $\phi\in\pi_2(\x,\x')$ with
$n_{\ws}(\phi)=0$ and
$2(n_{z^1_v}(\phi)+n_{z^2_v}(\phi))=-\epsilon_\x(v)+\epsilon_{\x'}(v)$.
In fact, this homology class is a union of disjoint bigons. In
particular, each Kauffman state $x$ is realized by a unique generator
$\x\in\Ta\cap\Tb$ which has maximal Alexander grading among all
generators realizing $x$, corresponding to the map $\epsilon(v)\equiv
+1$, for all $v$, cf.  Figure~\ref{fig:AlexanderMinimizers}.

Our aim in this section is to establish the following precise version
of Theorem~\ref{thm:States}:

\begin{theorem}
        \label{thm:MoreStates}
        Fix a decorated projection for a singular knot $K$ with the property
        that the distinguished edge $Q$ is just below a singular point.
        Let $C_d(K,s)$ be the free Abelian group generated
        by generalized Kauffman states with Alexander grading $s$
        and Maslov grading $d$. There is a differential
        $$\partial \colon C(K,s) \longrightarrow C(K,s)$$
        which carries $C_d(K,s)$ to $C_{d-1}(K,s)$, with
        $$H_d(C_*(K,s),\partial)\cong \HFb_d(K,s).$$
\end{theorem}

Before giving the proof, we establish some lemmas. Indeed, we find it
convenient to focus first on the case where $K$ is a planar
singular link (which is in fact the case where
Theorem~\ref{thm:MoreStates} is the most valuable).

\subsection{Planar singular links}

We analyze the Heegaard diagram in the case where $K$ is a planar singular
link (that is, admits an injective projection to the plane), establishing a
special case of Theorem~\ref{thm:MoreStates} in this case, from which
Theorem~\ref{thm:Planar} follows. First, we establish several lemmas.

\begin{lemma}
  \label{lemma:AlexanderGradings}
  Let $K$ be a planar singular link.  Let $x$ and $y$ be two
  generalized Kauffman states, and let $\x,\y\in\Ta\cap\Tb$ unique
  Alexander grading maximizing intersection points representing them.
  Then, the difference between the Alexander gradings of $\x$ and $\y$
  $($in the sense of Equation~\eqref{eq:AlexanderDifference}$)$
  coincides with the difference between the Alexander gradings $S(x)$
  and $S(y)$ of the Kauffman states $x$ and $y$.
\end{lemma}

\begin{proof}
  Given $x$ and $y$, it suffices to construct a particular homology
  class $\phi\in\pi_2(\x,\y)$ which satisfies
  \begin{equation}
    \label{eq:AlexanderDifferenceFormula}
    S(x)-S(y)=\sum_{i=1}^{\ell+1}(2n_{w_i}(\phi)-n_{z^1_i}(\phi)-n_{z^2_i}(\phi)).
  \end{equation}
  Indeed, it suffices to consider the case where there are no $D^-$
  corners in $x$ or $y$, in view of the following observation.
  Consider first the special case where $x$ and $y$ agree at all but
  $n$ corners, where $x$ is assigned $D^+$ and $y$ is assigned $D^-$.
  In this case, it is easy to find a locally supported homology class
  $\phi\in\pi_2(\x,\y)$ with $\sum n_{z_i}(\phi)= -n$, $\sum
  n_{w_i}(\phi)\equiv 0$, as illustrated in Figure~\ref{fig:DomPM}.

  \begin{figure}
    \mbox{\vbox{\epsfbox{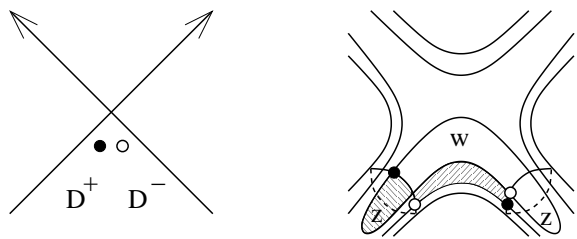}}}
    \caption{\label{fig:DomPM}
      {\bf{Domains for comparing $D^+$ and $D^-$.}}  We have
      illustrated on the right a domain $\phi$ which connects the
      Alexander grading maximizing generator of $D^+$ with that of $D^-$. All
      multiplicities are $+1$, $-1$, or $0$; regions with local
      multiplicity $+1$ are indicated by hatch marks from upper right
      to lower left, while regions with multiplicity $-1$ are
      indicated by the other hatching.}
  \end{figure}

  Suppose then that $x$ and $y$ contain no $D^-$ corners.  We
  construct a curve $\gamma$ which connects $\x$ to $\y$, as in
  Figure~\ref{fig:Domains}.  This curve $\gamma$ is a (possibly
  disconnected) closed path constructed from arcs within the $\alpha$
  and $\beta$-circles. In fact, if $\x$ has components $x_i$ and $\y$
  has components $y_i$, then $\gamma$ is constructed from arcs in
  $\alpha_i$ from $x_i$ to $y_i$ and arcs in $\beta_j$ from $y_j$ to
  $x_k$.
  
  Now, order the edges $\{e_i\}$ in the order they are encountered in
  the singular knot, starting from $Q$.  Let $v$ be some crossing
  where $e_i$ and $e_{i+1}$ meet.  Fix reference points $T_i$ on the
  Heegaard diagram, placed on top of the $i^{th}$ edge, see the lower
  left picture in Figure~\ref{fig:Domains}. We can draw an arc
  $\epsilon_i$ from $T_i$ to $T_{i+1}$, which crosses only one of the
  $\alpha$-circle which is the meridian for the $i^{th}$ edge, and
  none of the other circles. Similarly, let $\eta_i$ be the short arcs going
  from $w_v$ to $z^1_v$ and $z^2_v$ (inside the disk bounded by
  $\beta_v$). Then,
  \begin{align*}
    A(\x)-A(\y)&=\sum \#(\eta_i\cap \partial \phi)\\
    &= \sum\#((\epsilon_i+\eta_i)\cap \partial \phi),
  \end{align*}
  since $\cup_i \epsilon_i$ is a closed curve.
  We claim that
  \begin{align}
    \#(\epsilon_i+\eta_i)\cap \partial \phi &=
    S_v(x)-S_v(y)-\OneHalf\left(W_{x(v)}(K)-W_{y(v)}(K)\right),
    \label{eq:LocalContribution}
  \end{align}
  where here $W_r(K)$ denotes the winding number of $K$ around a point
  in the region $r$. We see this as follows.  The only part of
  $\partial\phi$ which intersects $\epsilon_i+\eta_i$ lies on the
  internal $\beta$-circle. Since the intersection number of
  $\epsilon_i+\eta_i$ with the internal $\beta$-circle is zero, it
  suffices to verify Equation~\eqref{eq:LocalContribution} replacing
  $\partial \phi$ with any arc connecting $\x$ to $\y$ on the internal
  $\beta$-circle. Thus, we can verify
  Equation~\eqref{eq:LocalContribution} by considering the various
  cases of $\x$ and $\y$ (locally, about each singular point), as
  illustrated in Figure~\ref{fig:Domains}.

  Now, summing Equation~\eqref{eq:LocalContribution}, we see that
  $$A(\x)-A(\y)=\left(S(x)-\OneHalf\sum_{v\in\Cr(P)} W_{x(v)}(K)\right)
  -\left(S(y)-\OneHalf\sum_{v\in\Cr(P)} W_{y(v)}(K)\right).$$
  But observe that 
  $$\sum_{v\in\Cr(P)} W_{x(v)}(K)=\sum_{r\in R(P)-X-Y} W_r(K)$$ is
  independent of the Kauffman state $\x$.  This finishes the proof.
  \begin{figure}
    \mbox{\vbox{\epsfbox{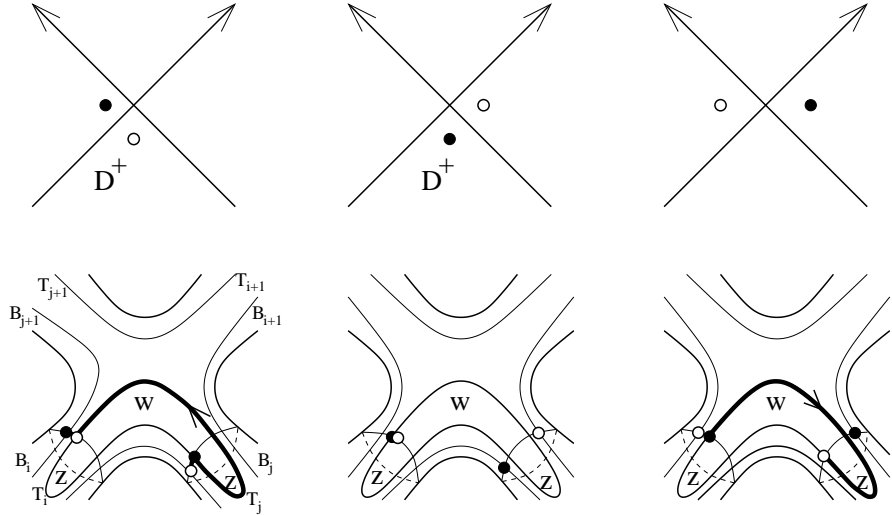}}}
    \caption{\label{fig:Domains} {\bf{Verifying
          Equation~\eqref{eq:LocalContribution}.}} We verify that
      equation, using arcs connecting the black dot (representing
      $\x$) to the white dot (representing $\y$) supported on the
      internal $\beta$-circle.}
  \end{figure}
\end{proof}

\begin{defn}
        \label{def:SingEquivalent} Let $K$ be a planar singular
        knot projection.  Two generalized Kauffman states $x$ and $y$
        are said to be {\em equivalent} if $x(v)\in \{A_v, D^-_v\}$ if
        and only if $y(v)\in\{A_v,D^-_v\}$.  Let $x$ be a generalized Kauffman
        state. This determines an associated subgraph of the knot
        projection as follows.  At each crossing $v$, if
        $x(v)\in\{A_v,D^-_v\}$, then we remove the lower left edge
        from the projection; if $x(v)\in\{C_v,D^+_v\}$, we remove the
        lower right edge from the crossing. We call the associated
        subgraph $\Gamma_x$ of the projection the {\em pruning}
        associated to the Kauffman state. Thus, two Kauffman states
        are equivalent if and only if they induce the same pruning.
\end{defn}

Note that the pruning is a graph, with exactly one in-coming edge at each
vertex, and at most  two out-going ones.

\begin{lemma}
        \label{lemma:ConnectedPruning}
        The pruning associated to a Kauffman state for a planar
        singular knot is a connected graph.
\end{lemma}

\begin{proof}
  Suppose that $\Gamma_x$ is not connected. We argue that $\Gamma_x$ must
  contain some cycle.  This is seen by taking some edge $e$ which is in a
  different path component from the initial point, and continuing backwards
  through $e$. Since each vertex in $\Gamma_x$ has a unique incoming edge, this
  process can be continued; so the component through $e$ must contain a cycle.
  
  Now, we consider some closed circuit $C$ in $\Gamma_x$ in a path
  component disjoint from the distinguished edge $Q$. The circuit $C$
  can be thought of as enclosing a region $\Delta$ of the projection
  which does not contain the distinguished edge. We restrict our
  Kauffman state to $\Delta$. The restriction of our projection to
  $\Delta$ gives the region, which is topologically a disk, the
  structure of a $CW$ complex, where the vertices correspond to the
  $V$ crossings, and edges corresponding to the $E$ arcs in $\Delta$,
  and which divide $\Delta$ into $F$ faces. There are four kinds of
  vertices in the one-complex of $\Delta$: bivalent ones (with one
  in-coming and one out-going edge) whose number is $D$, vertices with
  two in-coming and one out-going edge whose number is $T_1$, vertices
  with one in-coming and two out-going edges whose number is $T_2$,
  and four-valent ones (two in-coming, and two out-going edges) whose
  number is $W$, so that $V=D+T_1+T_2+W$.  Note that all vertices
  counted in $D$, $T_1$, and $T_2$ occur on the boundary of $\Delta$.
  By counting edges, we can verify that $T_1=T_2$. In fact, the total
  number of edges is given by $E=D+3T_2+2W$. But since $\Delta$ is a
  disk, its Euler characteristic is one, so we get the relation that
  $F-T_2-W=1$. Since each face in $\Delta$ is occupied by a Kauffman
  corner (we are using here the fact that $Q$ is not contained in
  $\Delta$), the restriction of our Kauffman state $x$ to $\Delta$
  demonstrates that $F\leq T_1+W$ (at each vertex $v$ counted in $D$
  or $T_2$, $x(v)$ is not a face of $\Delta$, whereas at each vertex
  $v$ in $W$, $x(v)$ is a face of $\Delta$, while at each vertex $v$
  in $T_1$, $x(v)$ might or might not be a face of $\Delta$),
  contradicting $T_1=T_2$ and $F-T_2-W=1$.
\end{proof}

This has the following easy consequence:

\begin{lemma}
  \label{lemma:EquivalentStates} 
If $x$ and $y$ are generalized Kauffman states for a planar singular
knot $K$ which are equivalent, then $x$ and $y$ coincide.
\end{lemma}

\begin{proof}
  We construct the following subset $P$ of the complement of
  $\Gamma_x$.  Suppose that $x$ and $y$ differ at a crossing $v$.
  Then, of course there is a different crossing $w$ with the property
  that $y(w)$ and $x(v)$ are assigned to the same region. We then
  connect $y(w)$ and $x(v)$ along some arc in the complement of the
  knot projection. Next, we connect $x(v)$ to $y(v)$ by an arc which
  crosses one of the four edges of our knot projection; but that is
  precisely the edge removed to obtain $\Gamma_x$. We continue this
  procedure.
  
  In this manner, we construct a collection of closed curves $P$. If
  $P$ is non-empty, let $R$ be any connected component of $P$. It is
  easy to see that $R$ divides the plane into two regions, both of
  which contain points in $\Gamma_x$. But this contradicts the fact
  that (by Lemma~\ref{lemma:ConnectedPruning}) $\Gamma_x$ is
  connected.
\end{proof}

We now have the ingredients required to establish
Theorem~\ref{thm:MoreStates}, at least in the case when $K$ is a
planar singular link.

\begin{prop}
        \label{prop:StateSums}
        Theorem~\ref{thm:MoreStates} holds for planar singular
        links; i.e. if $K$ is a planar singular link, let
        $C_d(K,s)$ be the free Abelian group generated by generalized
        Kauffman states with Alexander grading $s$ and Maslov grading
        $d$. There is a differential
        $$\partial \colon C(K,s) \longrightarrow C(K,s)$$
        which carries $C_d(K,s)$ to $C_{d-1}(K,s)$, with
        $$H_d(C_*(K,s),\partial)\cong \HFb_d(K,s).$$
\end{prop}

\begin{proof}
  If $\epsilon$ and $\epsilon'$ differ at a single singular point $v$,
  then the homology class $\phi$ clearly admits a single holomorphic
  representative up to reparametrization. In fact, by placing
  basepoints at all the other regions of our Heegaard diagram, we
  obtain a filtration of the chain complex $\CFb$ whose $E_0$ term
  counts only these short differentials; i.e. its differentials are
  given by tensor product of the space of Kauffman states with $\ell$
  chain compexes of the form
  $$\begin{CD} \Field [U_1,\ldots ,U_{\ell}]@>{U_i}>> \Field [U_1,
      \ldots ,U_{\ell}]
  \end{CD}.$$ It is easy to see that the $E_1$ term,
  now, is simply the free $\Field$-module generated by the Kauffman
  states.  Indeed, this homology is carried by the pairs
  $(x,\epsilon^-)$, where $\epsilon^-(v)=-1$ for each $v$. 
   
  It remains to identify the Alexander and Maslov gradings of the
  intersection points $(x,\epsilon^-)$ with the corresponding gradings
  for their underlying Kauffman states.  This statement for the
  Alexander grading is an immediate consequence of
  Lemma~\ref{lemma:AlexanderGradings}.
\begin{figure}
\mbox{\vbox{\epsfbox{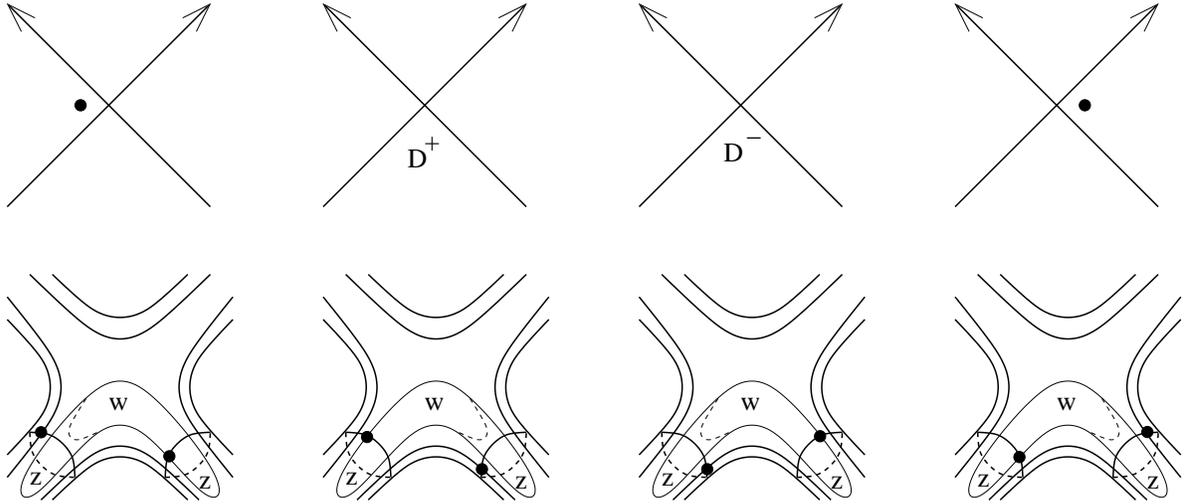}}}
\caption{\label{fig:IsotopeDiagram} {\bf{Isotopies of the standard
      diagram.}}  Given the generator $\x$, we consider the new
  diagram obtained by isotoping the central $\beta$-circle as
  indicated by the light dotted line. Note that the new diagram for
  $A$ and $D^-$ coincide, while the diagram for $C$ and $D^+$
  coincide.}
\end{figure}
For the statement about the Maslov gradings, it suffices to show that
the algebraic grading $N$ of any of the Alexander grading maximizing
generators vanishes identically.  To this end, let $\x$ be a generator
and $x$ its underlying Kauffman state. We claim that there is a new
Heegaard diagram, obtained by moving the central $\beta$-circles, each
across exactly one of the two $z^i_j$ for fixed $j$, so as to cancel
two of the intersection points of this $\beta$ with one of the two
$\alpha$-circle meridians, depending on the value of $x(v)$, as
illustrated in Figure~\ref{fig:IsotopeDiagram}. Let $\gamma_v$ denote
the isotoped image of $\beta_v$. Note that the $\gamma_v$ are chosen
so that the intersection point $\x$ persists into the new
diagram. Indeed, the generators for the new diagram now correspond to
pairs $(y,\epsilon)$, where here $y$ is some Kauffman state equivalent
to $x$, and $\epsilon\colon s(P)\longrightarrow \{\pm 1\}$. According
to Lemma~\ref{lemma:EquivalentStates} it follows at once that the
generators for the new diagram correspond simply to maps
$\epsilon\colon s(P)\longrightarrow \{\pm 1\}$, since they all have
the same underlying Kauffman state. We claim that this is an
admissible Heegaard diagram for $S^3$ with $\ell$ basepoints (supplied
by $\ws$). To see admissibility, we proceed as follows. Let $\Pi$
be a periodic domain with multiplicity $m$ near some singular point
$v$. Let $u$ be the preceding vertex to $v$ in the pruning
$\Gamma_x$. If $\Pi$ has only non-negative multiplicities, then its
local multiplicity near $u$ must be greater than or equal to
$m$. Iterating this, we can go back to the root of the pruning, in
view of Lemma~\ref{lemma:ConnectedPruning}, where the local
multiplicity in turn must equal zero. This shows that our periodic
domain must vanish identically.
  \begin{figure}
    \mbox{\vbox{\epsfbox{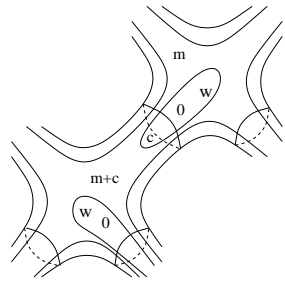}}}
    \caption{\label{fig:Admissibilty}
      {\bf{Verifying admissibility.}}  We have illustrated here the
      Heegaard diagram near two vertices, along with part of a
      periodic domain. (We have suppressed the basepoints $z$, as they
      no longer play a role.)  If the periodic domain has local
      multiplicity $m$ near the vertex $v$, and local multiplicity $c$
      in the small region, then its local multiplicity near the
      preceding vertex $u$ in the pruning must be $m+c$.}
  \end{figure}
  
  Next, note that the isotopy, of course, changes the Floer homology;
  but it leaves invariant the theory where all $U_i=1$ for $i>1$
  (inducing an isomorphism between $N$-graded theories).  Indeed,
  recall~\cite{HolDiskFour} that the isomorphism induced on homology
  by isotopies can be thought of as induced by a chain map which 
  counts holomorphic triangles. Moreover, the count of holomorphic
  triangles clearly takes the generator corresponding to
  $(x,\epsilon)$, to the generator corresponding to the same
  $\epsilon$ (respecting $N$-gradings). The result follows.
\end{proof}
This immediately gives a proof of Theorem~\ref{thm:Planar}:

\begin{proof}[of Theorem~\ref{thm:Planar}.]
  Applying Proposition~\ref{prop:StateSums}, we see that the
  differential actually must vanish identically, since the state sum
  ensures that for fixed Alexander grading, the generators of the
  chain complex have fixed Maslov grading.
\end{proof}

\subsection{The general case of Theorem~\ref{thm:MoreStates}.}

The discussion above can be generalized to the case of (non-planar) singular
knots, as well.  For example, we have the following generalization of
Lemma~\ref{lemma:AlexanderGradings}:

\begin{lemma}
  \label{lemma:AlexanderGradingsGen}
  Let $K$ be a singular link.  Let $x$ and $y$ be two generalized
  Kauffman states, and let $\x,\y\in\Ta\cap\Tb$ unique Alexander
  grading maximizing intersection points representing them.  Then, the
  difference between the Alexander gradings of $\x$ and $\y$ $($in the
  sense of Equation~\eqref{eq:AlexanderDifference}$)$ coincides with
  the difference between the Alexander gradings $S(x)$ and $S(y)$ of
  Kauffman states $x$ and $y$.
\end{lemma}

\begin{proof}
  We argue as in the proof of Lemma~\ref{lemma:AlexanderGradings},
  establishing Equation~\eqref{eq:LocalContribution} in the presence
  of non-singular intersection points.  This time, it is the central
  $\alpha$-circle which meets the curve $\cup\epsilon_i$. Its total
  intersection number with this curve is zero, so it sufficies to
  verify Equation~\eqref{eq:LocalContribution} by replacing
  $\partial\phi$ with arcs in the central $\alpha$-circle which
  connect the two generators $\x$ and $\y$, compare
  Figure~\ref{fig:HomotopyClass}. This is straightforward to verify,
  and the previous proof goes through.
\begin{figure}
  \mbox{\vbox{\epsfbox{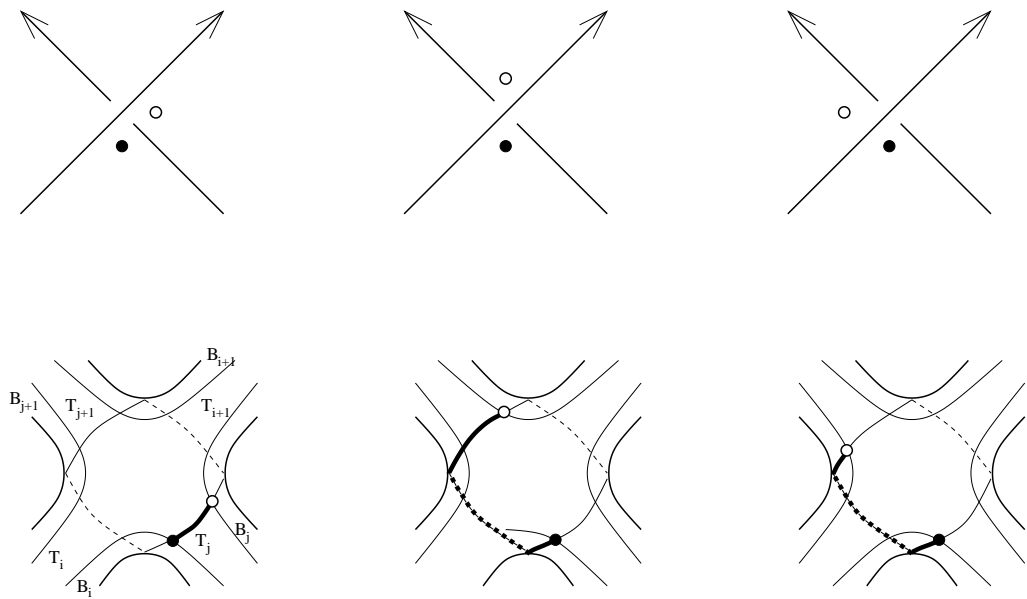}}}
\caption{\label{fig:HomotopyClass}
  {\bf{Verifying Equation~\eqref{eq:LocalContribution} in the presence
      of non-singular crossings.}} 
  The black dot represents $\x$, the white represents
  $\y$.  For negative crossings, rotate the
  diagram.}
\end{figure}
\end{proof}
\begin{proof}[of Theorem~\ref{thm:MoreStates}]
        Proceed as in the proof of Proposition~\ref{prop:StateSums}
        to see that there is a chain complex generated by Kauffman
        states.  Lemma~\ref{lemma:AlexanderGradingsGen} verifies the
        identification between the Alexander gradings coming from
        the Heegaard diagram with the state sum formula.
        
        It remains to verify that the Maslov grading for each
        Alexander grading maximizing generator $\x$ corresponding to a
        given Kauffman state $x$ is given by the state sum formula;
        indeed, since we have already verified the corresponding
        statement for the Alexander gradings, it suffices to verify
        the corresponding statement for the $N$-grading (which has the
        advantage that it is independent of the placement of the $z_i$
        basepoints).  Proceeding as before (as indicated in
        Figure~\ref{fig:IsotopeDiagram}), we isotope $\beta_v$ across
        one of the basepoints $z^i_j$ at each singular crossing $v$ to
        get a new circle. In fact, we can equivalently view this as an
        isotopy of one of the $\alpha$-circles intersecting $\beta_v$,
        replacing it with a new circle $\gamma_v$ disjoint from
        $\beta_v$.  Consider next a non-singular crossing $v$. If
        $x(v)$ is of type $A$ or $C$, we replace the $\alpha$-circle
        by a new circle $\gamma_v$ which is a meridian for the
        corresponding in-coming edge. If $x(v)$ is of type $D$, we
        replace its corresponding $\beta$-circle by the meridian
        $\gamma_v$ for either of the two in-coming edges.  Finally if
        $x(v)$ is of type $B$, we replace the $\alpha$-circle by a new
        circle $\gamma_v$ supported locally near the crossing
        pictured, as shown by Figure~\ref{fig:ThetaPrimed}.

\begin{figure}
\mbox{\vbox{\epsfbox{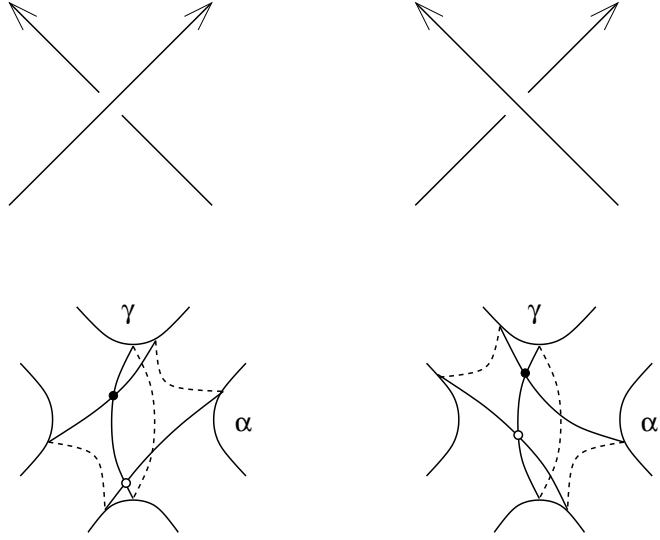}}}
\caption{\label{fig:ThetaPrimed} {\bf{A generator $\Theta'$.}}  The
  black dots represent the generator $\Theta'\in\Tc\cap\Ta$.  At each
  vertex, the $\gamma$- and $\alpha$-circles meet in two points, one of
  which we colored black, the other white. For a positive crossing,
  the black has $N$-grading one greater than the white generator; this
  can be verified using the obvious bigon (minus disks) going from
  black to white.  Similarly, at a negative crossing, the $N$-grading
  of the white generator is one greater than the $N$-grading of the
  black one.}
\end{figure}

        In this manner, we obtain a new diagram
        $(\Sigma,\gammas,\betas,\ws,\zs)$, equipped with a set
        $S=\Tc\cap\Tb$, from which we have a map to the Kauffman states
        for $K$.  Indeed, elements of $S$ map to the Kauffman states for
        the projection of another knot $K_0$ which is obtained by
        resolving all the crossings where $x(v)=B$. Let $x_0$ denote
        the induced Kauffman state on $K_0$. It is easy to see that
        any other intersection point  $\y$ of $\Tc\cap\Tb$ induces a
        Kauffman state $y$ which associates $B$ or $D$ to each vertex
        $v$ where $x(v)=B$, and hence it can be restricted to   $K_0$.
        In fact, $y_0$ is equivalent to $x_0$, and hence by
        Lemma~\ref{lemma:EquivalentStates} $x_0=y_0$. From this, it
        follows easily that $x=y$. Indeed,
        Lemma~\ref{lemma:ConnectedPruning} shows that the $\gamma$
        circles divide $\Sigma$ into $\ell+1$ components. Since the
        span of the $\gamma_i$ is contained in the span of the
        $\alpha_i$, it follows that these two spans coincide,
        and that $(\Sigma,\gammas,\betas,\ws)$ is a Heegaard
        diagram for $S^3$.
        
        Now, given generator $\x\in\Ta\cap\Tb$ (maximizing Alexander
        grading among all intersection points corresponding to the
        Kauffman state $x$), we have exhibited a new Heegaard diagram,
        equipped with a corresponding Alexander grading maximizing
        generator $\x'\in\Ta\cap\Tc$, which clearly has $N$-grading
        equal to zero. In fact, $\Ta\cap\Tc$ also has a minimal number
        of intersection points. Let $\Theta$ denote its $N$-maximizing
        generator. If $n$ denotes the number of negative crossings $v$
        where $x(v)=B$ (and $x$ is the Kauffman state corresponding to
        $\x$), and $p$ denotes the number of such positive crossings,
        then we claim that there is an obvious
        $\psi\in\pi_2(\Theta',\x,\x')$ where
        $\gr(\Theta')=\gr(\Theta)-n$, and $\Mas(\psi)=-p$,
        cf. Figure~\ref{fig:ThetaTriangle}.  The map induced by
        counting holomorphic triangles preserves gradings, in the
        sense that $N(\Theta')+N(\x)-\Mas(\psi)=N(\x')$, from which it
        follows now that $N(\x)=n-p$. Comparing with
        Figures~\ref{fig:MasGradings} and~\ref{fig:AlexGradings}, we
        have verified that $N(\x)=M(x)-2S(x)$.
\begin{figure}
\mbox{\vbox{\epsfbox{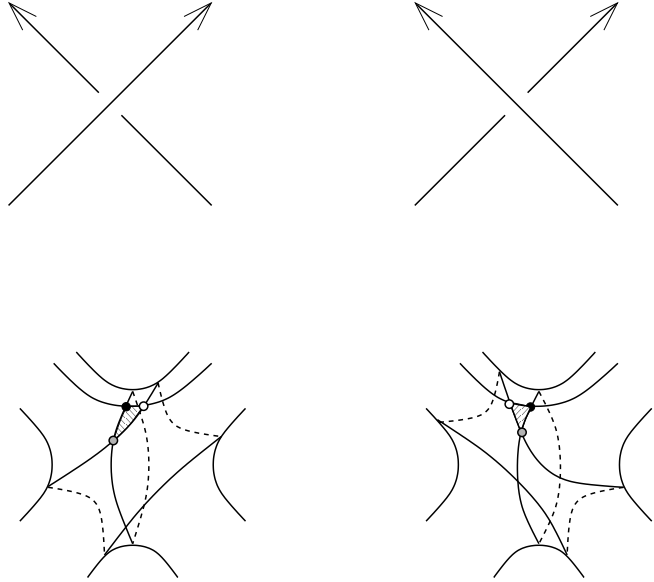}}}
\caption{\label{fig:ThetaTriangle} {\bf{Triangles.}}  We denote the
  generator $\Theta'$ by a grey dot;  $\x'\in\Ta\cap\Tc$ and
  $\x\in\Ta\cap\Tb$. The triangle on the left (corresponding to the
  positive crossing) has local coefficent $-1$ (and indeed Maslov
  index equal to $-1$) while the one on the right has positive local
  multiplicity (and Maslov index equal to zero).}
\end{figure}
\end{proof}
The proof of Theorem~\ref{thm:Euler} is an immediate corollary of
Theorem~\ref{thm:MoreStates}:

\begin{proof}[of Theorem~\ref{thm:Euler}.]
  Equation~\eqref{eq:EulerHFb} is an immediate consequence of
  Theorem~\ref{thm:MoreStates} and Proposition~\ref{prop:StateSumFormula}.

  Equation~\eqref{eq:EulerHFa} follows by a simple comparison of the
  chain complexes: each generator $\x$ for $\CFa$ corresponds
  infinitely many generators for $U_1^{n_1}\cdot...\cdot
  U_\ell^{n_\ell}$ indexed by $\ell$-tuples of non-negative integers
  $(n_1,...,n_\ell)$.  Each such generator occupies Alexander grading
  $A(\x)-\sum_{k=1}^{\ell} n_k$ and $N$-grading equal to the
  $N$-grading of $\x$. Equation~\eqref{eq:EulerHFa} follows.
\end{proof}

According to Theorem~\ref{thm:Planar}, calculating $\HFb(K)$ for a
planar knot is equivalent to computing its Alexander polynomial. 
This can be efficiently done either using the state sum formula
$$\Delta_{K}(T)=\sum_{x\in K(P)} \prod_{v\in \Cr(P))} Z_v(x)$$
we have discussed in Section~\ref{sec:Alex}, or the
skein relations
\begin{eqnarray*}
      \Delta_{K^+}(T)&=&\Delta_{K} (T)+T^{\OneHalf}\cm {\Delta}_{K^o}(T)
       \\ \Delta_{K^-}(T)&=& \Delta_{K}(T)+T^{-\OneHalf}\cm\Delta_{K^o} (T), \\
\end{eqnarray*}
we have used in the definition of the Alexander polynomial for singular
links. 

\section{Some calculations}
\label{sec:Planar}

\subsection{Planar singular links}
\begin{figure}
\mbox{\vbox{\epsfbox{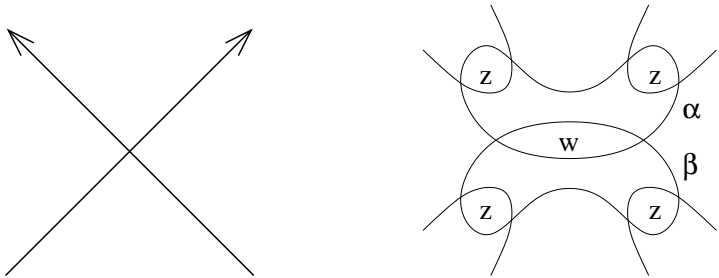}}}
\caption{\label{fig:Plancross} {\bf{Heegaard diagram of a planar singular
      link.}} The crossing on the left (obtained by contracting a thick edge)
	is replaced by the piece of Heegaard diagram on the right.}
\end{figure}

Recall that a singular link $K$ is called \emph{planar} if it admits an
injective projection to the plane. For such links a simple Heegaard diagram of
genus zero can be given in the following way.  Fix a planar singular link $K$,
consider an injective projection, contract all its thick edges to singular
points and take the $\alpha$- and $\beta$-curves at every crossing as it is
instructed by Figure~\ref{fig:Plancross}.  Note that an $\alpha$-curve
corresponds to each thick edge and the two outgoing thin edges, while a
$\beta$-curve corresponds to a thick edge and two incoming thin edges.
By taking these curves for all the crossings and then deleting an (arbitrary)
$\alpha$- and $\beta$-curve from the collection, we get a Heegaard diagram of
a singular link.  

\begin{prop}
The resulting Heegaard diagram is compatible with the given planar
singular link $K$. \qed
\end{prop}

We can distinguish coordinates of an intersection point $\x =(x_1, \ldots ,
x_{\ell })\in \Ta\cap \Tb$ according to whether $x_i$ corresponds to a thin or
thick edge --- in the diagram it is reflected by the fact whether $x_i$ is
near a base point of type $\z$ or of type $\w$.  Since near a fixed $w$ or $z$
the $\alpha$-- and $\beta$--curves intersect each other in two points, we can
group our intersection points into groups of cardinality $2^{\ell}$.

\subsection{An example}
We will illustrate the above principle by an example.  This example
also shows that the Floer homology theory $\HFa$ is not determined
naively by the Alexander polynomial of a planar singular knot.

\begin{figure}
\mbox{\vbox{\epsfbox{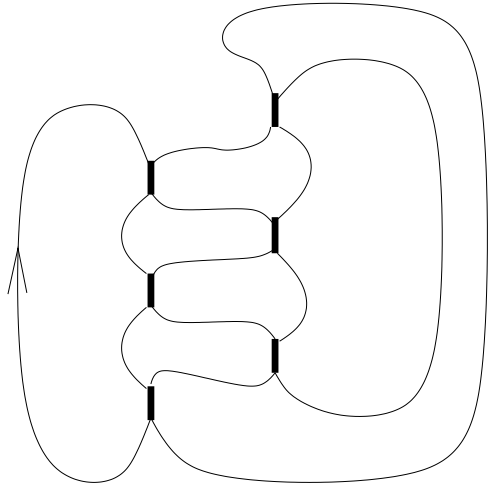}}}
\caption{\label{fig:Toruslink}
{\bf{Singular knot given as the singularization of the projection of the (3,3)
    torus link.}}}
\end{figure}

\begin{figure}
\mbox{\vbox{\epsfbox{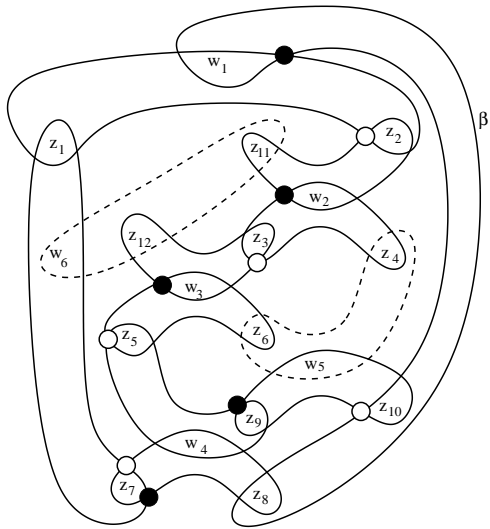}}}
\caption{\label{fig:PlanarDiag}
{\bf{The planar Heegaard diagram corresponding to the singular knot
projection of Figure~\ref{fig:Toruslink}.}}}
\end{figure}

By taking the (3,3) torus link and singularizing its natural
projection we get the singular knot $K$ depicted by Figure~\ref{fig:Toruslink}.
The planar Heegaard diagram corresponding to this projection is shown by
Figure~\ref{fig:PlanarDiag}.  The diagram also indicates (with dashed lines)
the $\alpha$- and $\beta$-curves which we delete according to the algorithm
for constructing the diagram from the projection.

The Alexander polynomial of the singular knot $K$ given above can
easily be computed from the state sum formula, giving 
\[
\Delta _K (T)=T^2+5 \cm T+9+5\cm T^{-1}+T^{-2}.
\]
It follows from Theorem~\ref{thm:Euler}
that
\[
\sum_{s,d}(-1)^d \rk~\HFa_d(K,s) \cm T^s =
-T^7+6 \cm T^5-21 \cm T^3+21 \cm T^2-6+T^{-2}
\]

In the following we will identify the generators of the Floer chain
complex $\CFa$ and determine part of the boundary map $\Da$ for this
particular singular knot. In certain Alexander gradings, this 
differential is trivial; indeed, we have the following:

\begin{prop}
  The rank of $\HFa(K, s)$ for $s=2,3$ is $21$, for $s=0,5$, it is $6$, 
  for $s=-1,6$ is $0$ while for $s=-2,7$ it is $1$.
\end{prop}

\begin{proof}
This is a straightforward calculation, explicitly identifying the
generators, and finding some boundary maps. See the proof of
Proposition~\ref{prop:EulerZero} for more details.
\end{proof}

The rank of $\HFa$ is not equal to its Euler characteristic.  More
precisely we show that for the above singular link $K$

\begin{prop}
  \label{prop:EulerZero}
The Floer homology group $\HFa_*(K, 4)$ is nontrivial, although its
Euler characteristic is zero.
\end{prop}

\begin{proof}
In order to prove this statement, we first list the intersection points of the
tori $\Ta$ and $\Tb$ in $Sym^5(S^2)$ and identify the ones which have
Alexander grading $4$. First of all notice that there are four
types of intersection points: the points of type $A$ consist of those
intersection points which have coordinates near $(w_2,w_3, w_4, z_1, z_{10})$;
for type $B$ points this vector is $(z_2,z_3,z_5,z_7,z_{10})$, for type $C$
points this vector is $(w_2,w_3,z_1,z_8,z_9)$, and finally for type $D$ points
this vector is $(w_1,w_2,w_3,z_7,z_9)$.  The two possible coordinates near a
base point $z_i$ (or $w_j$) will be distinguished based on the property
whether it is to the left or to the right from the edge corresponding to the
base point (when using the orientation opposite of the vertex). Near a 
$w$-type point the two choices will be denoted by $L$ and $R$, near a $z$-type
point by $l$ and $r$.  We also keep the order of listing $w$'s first (with
increasing indices) followed by $z$'s (also with increasing indices).  For
example, the intersection point $\x _1 \in \Ta\cap \Tb$ indicated by the heavy
dots in Figure~\ref{fig:PlanarDiag} is represented by $D(R,L,L,r,l)$, while
the light dots $\y\in\Ta\cap\Tb$ in the same figure is represented by  $B(l,r,l,l,l)$.

It is easy to see that changing $L$ to $R$ drops the $N$-grading by
$1$ and raises the Alexander grading by $2$, while a change of $l$ to
$r$ drops the $N$-grading by $1$ and raises the Alexander grading by
$1$.  For example, the highest Alexander grading is attained by the
intersection point $D(L,L,L,l,l)$ while the lowest Alexander grading
is taken by $A(R,R,R,r,r)$.

The above combinatorial
computation shows that there are twenty intersection points with
Alexander grading $4$, ten of which have (up to a suitable
translation) Maslov grading $2$ and ten have Maslov grading $1$.

Consider the following eleven intersection points: $\x_1=D(R,L,L,r,l)$
and $\x_2 = B(l,r,l,l,l)$ (depicted in Figure~\ref{fig:PlanarDiag} by
the heavy and light circles), and also $\b_1=B(r,l,l,l,l)$,
$\b_3=B(l,l,r,l,l)$, $\b_4=B(l,l,l,r,l)$, $\b_5=B(l,l,l,l,r)$, and
$\a_3=A(L,L,R,l,l)$ (of Maslov grading 1); and $\D_2 =D(L,R,L,r,l),
\D_3 =D(L,L,R,r,l), \D_5=D(L,R,L,l,r)$, and $\D_6 = D(L,L,R,l,r)$ (of
Maslov grading 2). (In this notation $\x_1=\D_1$, and $\x_2=\b_2$.)
There are nine nonnegative homotopy classes with $n_{\zs}=n_{\ws}=0$
connecting nine of the above eleven intersection points $\b_i$ ($i
\neq 2$), $\D_j $ ($j \neq 1$) and $\a_3$ to the remaining nine
intersection points of Alexander grading $4$. It is easy to see that
the contribution of each of these nine homotopy classes in the
boundary map is equal to $\pm 1$.  Next we will show that there are no
more pairs of intersection points among the eleven above with the
property that a homotopy class connecting them with
$n_{\zs}=n_{\ws}=0$ exists. Clearly, this fact implies that there are
no more boundary maps we should take into account when computing $\HFa
(K, 4)$, verifying that $\HFa (K, 4)=\Field \oplus \Field$ (with
$\Field = \Z /2\Z$), concluding our computation.

\begin{figure}
\mbox{\vbox{\epsfbox{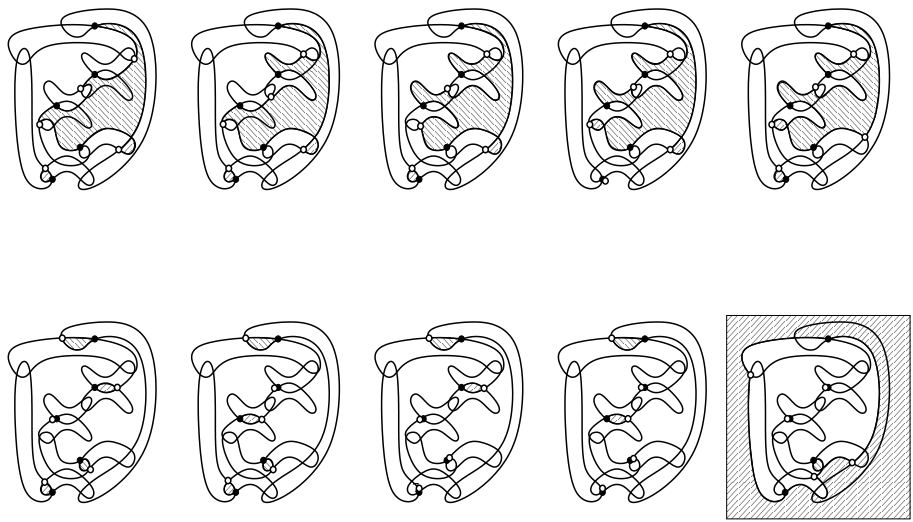}}}
\caption{\label{fig:FindDomains} {\bf{Domains}.}  We depict domains
  from $\x_1$ to the other $10$ distinguished generators (in order:
  $\b_1$, $\x_2$, $\b_3$, $\b_4$, $\b_5$, $\D_5$, $\D_6$, $\D_2$,
  $\D_3$, $\a_3$).  The generator $\x_1$ is denoted by the black
  dots. Local multiplicities are $0$, $-1$ (denoted by hatchings in
  one direction), and $+1$ (hatchings in the other direction).}
\end{figure}

For the last claim about the nonexistence of further homotopy classes with
$n_{\zs}=n_{\ws}=0$ we argue as follows. We find domains from $\x_1$ to the
other ten distinguished intersection points
$\{\b_1,\x_2,\b_3,\b_4,\b_5,\a_3,\D_2,\D_3,\D_5,\D_6\}$ (for example, the ones
pictured in Figure~\ref{fig:FindDomains}).  Associate to each such domain
$\cald $ the corresponding 18--dimensional vector $(n_{\zs}(\cald ),
n_{\ws}(\cald) )$.  (Recall that there are 12 points of type $z$ and 6 of type
$w$.)  This 18-dimensional vector space contains a subspace $V$ generated by the
vectors $(n_{\zs}({\mathcal P}),n_{\ws}({\mathcal P}))$, where ${\mathcal P}$
runs over all domains with boundary among the $\alpha_i$ and $\beta_j$. Our
claim amounts to showing that the original $10$ vectors, together with
the additional zero vector, are distinct modulo $V$.  To this end, consider
functions $F_1=w_1^*-z_1^*-z_8^*+z_7^*,$ $F_2=w_2^*-z_3^*-z_{11}^*+z_{12}^*$, 
$F_3=w_3^*-z_6^*-z_3^*+z_4^*$, $F_4=w_4^*-z_8^*-z_9^*+z_{10}^*$,
$F_5=w_5^*-z_6^*-z_9^*+z_5^*$, which can be easily shown to vanish on $V$.  It
is straightforward to verify that these five functions take on $11$ distinct
values on the $11$ vectors.  (For the sake of completeness, we list the
vectors of our chosen domains connecting $\x_1$ to the points in the obvious
basis of $\Z ^{18}$, where the basis vectors are identified with the $z_i$'s
and the $w_j$'s: for $\b_1$ it is $(z_4+z_6-z_7-z_{10})$, $\b_2$ is
$(z_6-z_7-z_{10}+z_{11})$, $\b_3$ is $(-z_7-z_{10}+z_{11}+z_{12})$, $\b_4$ is
$(-z_5-z_{10}+z_{11}+z_{12})$, $\b_5=-z_5-z_7+z_{11}+z_{12}$, $\D_5$ is
$-z_7+z_9+w_1-w_2$, $\D_6$ is $-z_7+z_9+w_1-w_3$, $\D_2$ is $w_1-w_2$, $\D_3$
is $w_1-w_3$, $\a_3$ is $z_9-z_{10}$.)

We conclude that $\HFa(K,4)$ is generated by the two generators
$\x_1=D(R,L,L,r,l)$ and $\x_2 = B(l,r,l,l,l)$ depicted in
Figure~\ref{fig:PlanarDiag} by the heavy and light
circles. Consequently we see that
\[
\HFa (K, 4)\cong \Field \oplus \Field ,
\] 
showing that the total rank of the Floer homology of $\HFa$ can exceed
the absolute value if its Euler characteristic.
\end{proof}

Note that a similar
calculation can be performed for $s=1$, showing that the corresponding
Floer homology group is also isomorphic to $\Field \oplus\Field$.

\commentable{ }
\bibliographystyle{plain} \bibliography{biblio} 

\end{document}